\documentclass[a4paper,11pt,nonatbib]{article}

\usepackage[utf8]{inputenc}
\usepackage[backend=biber,style=apa]{biblatex}
\usepackage{amsmath, amssymb, amsfonts}
\usepackage{geometry}
\usepackage{mathtools}
\usepackage{bm}
\usepackage{graphicx} 
\usepackage{hyperref}
\usepackage{xspace}
\usepackage{xcolor}
\usepackage{color} 
\usepackage{colortbl}
\usepackage{todonotes}
\usepackage{makecell}
\usepackage{standalone}
\usepackage{subcaption}
\usepackage{glossaries}
\usepackage{rotating}

\usepackage[ruled,vlined,linesnumbered,noresetcount]{algorithm2e}
\RestyleAlgo{ruled}
\SetKwComment{Comment}{/* }{ */}

\usepackage{authblk}
\usepackage{booktabs}
\usepackage{multirow}
\usepackage{comment}
\usepackage{csvsimple}
\usepackage{ifthen}
\usepackage{siunitx}
\usepackage[normalem]{ulem}
\usepackage{ulem}
\usepackage{pgfplots}
\usepackage{pgfplotstable}
\usepgfplotslibrary{groupplots}
\pgfplotsset{compat=1.18}

\definecolor{colBL}{HTML}{4878d0}
\definecolor{colMcCormick}{HTML}{ee854a}

\newcommand{\nlmodel}{\textsc{NL}\xspace} 
\newcommand{\pwlmodel}{\textsc{PWL}\xspace} 
\newcommand{\pwlBLmodel}{\textsc{BL}\xspace}
\newcommand{\heurShort}{\textsc{Heur}\xspace} 

\newcommand{\equidist}{\textsc{equidist}\xspace} 
\newcommand{\opt}{\textsc{opt}\xspace} 

\newcommand{\TPCodeFromInstance}[1]{%
\ifnum\pdfstrcmp{#1}{Athens}=0
ATH%
\else\ifnum\pdfstrcmp{#1}{Grid}=0
GRD%
\else\ifnum\pdfstrcmp{#1}{Lowersaxony}=0
LOS%
\else\ifnum\pdfstrcmp{#1}{Sioux}=0
SIO%
\else
UNK%
\fi\fi\fi\fi
}

\newacronym{MCF}{MCF}{Multi-Commodity Flow}
\newacronym{SAMCF}{SAMCF}{Service-Aware Multi-Commodity Flow}
\newacronym{QoS}{QoS}{Quality of Service}
\newacronym{NLP}{NLP}{Nonlinear Program}
\newacronym{MILP}{MILP}{Mixed Integer Linear Programming}
\newacronym{TAP}{TAP}{Traffic Assignment Problem}

\title{A Column Generation-based Fixed-point Heuristic for the Service-Aware Multi-Commodity Flow Problem}

\author[1]{Siv Marie Cartland Hansen}
\author[2]{Richard Martin Lusby}
\affil[1]{Technical University of Denmark \\ \texttt{scaha@dtu.dk}}
\affil[2]{Technical University of Denmark \\ \texttt{rmlu@dtu.dk}}

\date{}

\begin{document}
	
	\maketitle
        
        \begin{abstract}
    We study the Service-Aware Multi-Commodity Flow (SAMCF) problem, in which demand is elastic and governed by a logit choice model while routing is subject to hard capacity constraints. In a centralized, system-optimal setting, the network operator jointly determines how much demand to serve and how to route it. We formulate the SAMCF as a nonlinear program and propose an iterative fixed-point heuristic that alternates between solving an inelastic MCF via column generation and updating demand from the resulting service levels. Two linear approximations based on piecewise-linear demand functions and McCormick envelopes serve as benchmarks,  while a piecewise-linear outer-approximation of the demand function is used to provide valid lower bounds. Computational experiments on public transport instances show that the heuristic finds near-optimal solutions in under two seconds---orders of magnitude faster than the benchmark methods---while matching their solution quality on all instances they can solve within a ten-minute time limit.\\

    \end{abstract}

    \section{Introduction}

The \gls{MCF} problem is well‑studied in the Operations Research literature with diverse applications including telecommunications, transportation, and logistics, see~\textcite{salimifard2022a}. Given a capacitated, directed network consisting of a set of nodes $N$ and set of arcs $A$, the problem seeks to simultaneously route a set of commodities $D$ through the network. Each commodity $d\in D$ has an origin node, a destination node, and a demand volume $w_d$, while all commodities share the same underlying network capacities. The central challenge is to find feasible routing decisions that respect arc capacity constraints and which optimize a given objective. In the standard node-arc formulation, this amounts to routing each commodity $d\in D$ so as to minimize total cost $\sum_{d\in D}w_d\sum_{a\in A}t_af_{da}$, subject to flow conservation at every node and the joint arc capacity constraints $\sum_{d\in D}w_df_{da}\le c_a$ for each arc $a\in A$. Here $f_{da}\ge 0$ denotes the fractional flow of commodity $d$ on arc $a$, $t_a$ denotes the unit cost of using arc $a$, and $c_a$ denotes the capacity of arc $a$.

In many real-world applications, demand is inherently responsive to the level of service provided by the network. In transportation systems, travelers may choose between public transport and alternative modes based on travel time or cost, while in logistics networks, a customer's willingness to order may depend on delivery performance. This motivates the need for {\em service-aware} models that explicitly capture the interaction between network performance and demand. Introducing service-awareness fundamentally changes the structure of the \gls{MCF} problem, as routing influences service quality, which in turn affects the amount of demand attracted to the network. This interdependence between supply and demand leads to a coupled optimization problem that may be nonlinear.

In this paper, we study the \gls{SAMCF} in which demand is elastic and governed by a choice model, while routing decisions are subject to hard capacity constraints. Demand is elastic in the sense that a portion of the demand associated with a commodity may choose an outside alternative when network performance becomes unattractive. Unlike classical elastic-demand or user-equilibrium models, we consider a centralized, system-optimal setting in which the network operator determines both how much demand to serve and how to allocate capacity across commodities. Importantly, not all attracted demand must be accommodated, reflecting realistic capacity limitations that are typically found in practice. 
We model the \gls{SAMCF} problem as an \gls{NLP} and design an iterative fixed-point heuristic to solve it. The proposed heuristic alternates between solving an inelastic \gls{MCF}, in which demand is fixed, using column generation, and updating demand based on the experienced service levels. Through extensive computational experiments on diverse data sets, we demonstrate that the proposed heuristic is highly effective, consistently finding near-optimal solutions in seconds and significantly outperforming exact methods applied to the full formulation or its approximation. 

The contributions of the paper are fourfold. From a methodological perspective, the primary contribution is the development of an extremely competitive heuristic for the \gls{SAMCF} problem. In addition, we provide a comprehensive comparison with three exact methods for solving the nonlinear formulation, including the full formulation, a partially linearized approximation, and a fully linearized approximation. Not only is the proposed heuristic orders of magnitude faster, but it also provides solutions of --- in some cases substantially --- better quality as the benchmark approaches time out. In order to assess the quality of the heuristic solutions, we investigate a piecewise-linear outer approximation of the demand function and demonstrate that it provides substantially stronger lower bounds than those obtained by directly solving the full nonlinear formulation.  Conceptually, the notion of ``service-awareness'' has appeared in various contexts with different interpretations. We therefore provide a broad overview of the related literature and position our work within these streams. Finally, through our computational study, we offer perspectives on the role of demand elasticity in capacity-constrained networks. 

The structure of the paper is as follows. Section~\ref{sec:literature} summarizes studies that address service-awareness in the context of routing-related applications. Section~\ref{sec:problem} formulates the \gls{SAMCF} as a nonlinear optimization problem and provides an illustrative example. Section~\ref{sec:methods} presents different methodologies, including the proposed iterative fixed-point heuristic, a linear approximation of the proposed model in which the demand function is approximated by a piecewise-linear function, and a fully linear model in which the supply function is additionally approximated by McCormick envelopes. The performance of all approaches is then analyzed and compared in Section~\ref{sec:results} using an application in public transportation. Conclusions are given in Section~\ref{sec:conclusion}.

\section{Literature}
\label{sec:literature}

Recent research highlights the importance of modeling supply-dependent (endogenous) demand in network 
systems, where service decisions influence the demand that can be attracted and served across applications such as transportation, logistics, and telecommunications. \textcite{frejingerPerspectivesOptimizingTransport2025} emphasize the need for optimization frameworks that account for these demand--supply interactions, noting the limitations of models with exogenous demand. In this paper, we integrate demand response into capacitated \gls{MCF} models, viewing the \gls{SAMCF} as a network flow-based formulation within this broader line of research.

A growing body of literature has addressed supply-dependent demand across different application domains. These works differ primarily in how the interaction between service supply and demand is modeled, ranging from adoption and routing behavior in transportation planning, to demand elasticity in logistics, and service-level constraints in telecommunications.

In transportation planning, supply-dependent demand is often studied in the context of transit network design, where passenger demand responds to service quality through routing or adoption decisions. 
\textcite{tsaggourisQoSawareMulticommodityFlows2006} study a {\em \gls{QoS}-aware} variant of the \gls{MCF} problem motivated by transportation applications, where demand is sensitive to service quality. They model demand responsiveness through a consumption function that scales flow based on path-level deviations from shortest travel times, thereby embedding a form of elasticity directly into the flow formulation. 
This formulation relies on path-based flow representations, with demand responsiveness defined
at the path level to ensure decomposability, and is therefore not directly applicable to settings
 where elasticity is defined in terms of aggregated service measures. The problem is solved using a Lagrangian relaxation with column generation, where capacity constraints are dualized, and the pricing problem involves non-additive shortest paths that are solved approximately using a fully polynomial-time approximation scheme. A feasible primal solution is obtained by recovering flows from the approximately optimal dual solution, though capacity is not enforced explicitly and feasibility is achieved only a posteriori. Among existing approaches, this work is closely related to ours in that it embeds service quality directly within a \gls{MCF} framework. However, our approach differs in several key aspects: we do not require decomposability of the demand model, we maintain primal feasibility with respect to capacity constraints throughout, and we accommodate demand response through a logit-based model defined on aggregate service levels.

\textcite{klierUrbanPublicTransit2015} consider public transit network design with elastic demand and present an integer program that maximizes the expected number of public transit passengers. To endogenize passenger uptake, their response is governed by a binary logit model, where the public transport choice probability for a passenger depends solely on the best-established transit path. Capacity constraints are ignored, and the approach relies on an a priori path-enumeration step to identify reasonable paths for passengers.  \textcite{de-los-santosRailwayLineFrequency2017} integrate mode choice into a railway line frequency and train capacity optimization problem, using a linearized logit model to capture the modal split between the railway network and a competing mode representing alternative transport options. The authors assume that all passengers between the same origin-destination (OD) pair use the same path. A scalable approach for maximizing ridership in transit network design is proposed by~\textcite{bertsimasDataDrivenTransitNetwork2021}, where ridership depends on the number of opened lines, 
their operating frequencies, and the necessity of transferring. \textcite{cancaGeneralRapidNetwork2016,cancaAdaptiveNeighborhoodSearch2017,cancaIntegratedRailwayRapid2019} model demand elasticity using logit-based probabilities within an integrated network design and line planning problem, where infrastructure, frequencies, and demand are jointly determined. \textcite{cancaGeneralRapidNetwork2016} allows flows across different paths but elasticity is still largely driven by the best path.
\textcite{cancaAdaptiveNeighborhoodSearch2017} enforce exclusive shortest-path assignment, and  \textcite{cancaIntegratedRailwayRapid2019} consider a multi-path assignment, where elasticity depends on how an origin-destination pair is routed across the whole network, and present a two-level local search matheuristic that is able to solve instances of realistic size.

More recently, \textcite{hartlebModelingSolvingLine2023} have proposed a \gls{MILP} model for railway line planning that integrates mode and route choice using a logit-based demand model. The approach estimates passenger demand endogenously from traveler preferences while accounting for capacity constraints, on a precomputed set of feasible passenger paths. The resulting formulation is linearized by restricting choices to a discrete path set and representing service quality via the number of available routes, allowing the nonlinear logit-based demand response to be precomputed and embedded in a \gls{MILP} model. 

More recent approaches explicitly formulate and solve bilevel optimization problems that couple network design decisions with routing and behavioral responses. Such problems often exhibit a bilevel structure, where the upper level determines the network design and the lower level captures routing and/or demand outcomes. For example, \textcite{guanPathBasedFormulationsDesign2024,guanBilevelOptimizationHeuristic2025} and \textcite{basciftciCapturingTravelMode2023} consider an on-demand multimodal transit system design problem with endogenous adoption behavior, where rider decisions depend on the service quality induced by the network. The three papers essentially differ in their solution approaches. \textcite{guanBilevelOptimizationHeuristic2025} present different heuristics that iterate between network design and adoption decisions, whereas  \textcite{basciftciCapturingTravelMode2023} detail an exact decomposition of the problem using Benders decomposition. \textcite{guanPathBasedFormulationsDesign2024} instead propose a path-based reformulation that embeds routing and adoption decisions in a single-level model, enabling improved scalability to large instances. \textcite{weiTransitPlanningOptimization2022} consider a closely related setting by optimizing public transit schedules under competition from ride-hailing services and the presence of traffic congestion. Passenger behavior is modeled through a discrete choice framework, where mode choice depends on service attributes such as travel time and cost. Their model captures the interaction between service design, congestion, and passenger mode choice through a mixed-integer nonlinear formulation. The problem is solved using a bilevel heuristic algorithm that alternates between an upper-level transit scheduling problem and a lower-level user-equilibrium traffic assignment model. Bilevel approaches emphasize decentralized user responses and adoption decisions, often relying on pre-generated path sets or equilibrium-type assumptions, rather than explicitly controlling how much demand is served within the network.

In logistics and distribution systems, the interaction between supply and demand is often modeled through explicit demand elasticity, where demand depends endogenously on service attributes such as delivery time, reliability, or network coverage. For instance, \textcite{liuIterativePredictionandOptimizationELogistics2022} propose an iterative predict-and-optimize framework in which a machine learning model estimates demand as a function of service characteristics and optimization decisions are updated accordingly. Such approaches capture how improved service levels attract additional demand, typically through iterative schemes linking demand prediction and network optimization. Similarly, \textcite{zetinaProfitorientedFixedchargeNetwork2019} study a fixed-charge multicommodity network design problem in which demand is modeled via a gravity function depending on shortest-path costs in the induced network. Their model assumes an uncapacitated setting and abstracts from flow-level routing decisions, focusing instead on path selection and demand estimation. It is solved using an iterative matheuristic that combines column generation and slope scaling, which is shown to outperform directly solving 
the full model. \textcite{martin2021a} also study a profit‑maximizing service network design problem with endogenous demand, where delivery-time offerings and prices jointly determine shipment volumes and operations. It proposes an integrated model combining customer choice, pricing, and network design, with pre-enumerated path sets and a decomposition-based heuristic linking demand generation to routing decisions. Elasticity is captured via utility-based thresholds, where customers switch discretely between services or opt out as delivery time-price trade-offs change, yielding piecewise demand responses. \textcite{bilegan2022a} consider a  problem in service network design with endogenous demand in a simpler revenue-management setting for barges, where service levels are represented by delivery-time thresholds (due dates) and demand is determined through acceptance and routing decisions subject to these constraints; the problem is formulated as a \gls{MILP} model that jointly selects services, routes flows, and decides which demands to accept (fully or partially), and is solved directly using a commercial solver. Compared to transportation planning models, works in logistics place greater emphasis on demand estimation and forecasting, while often abstracting from detailed flow-level routing decisions and capacity allocation across commodities. 

In telecommunications and network routing, service-awareness is typically captured through \gls{QoS} requirements rather than explicit demand elasticity. In this setting, demand is typically given, but only flows satisfying certain service thresholds, such as, e.g., latency, bandwidth, or path length, are admitted or prioritized, see e.g.,~\textcite{khanGuaranteeingEndtoendQoS2021}. For example, \textcite{seliuchenkoEnhancedMulticommodityFlow2016} consider Software-Defined Networks and extend \gls{MCF} models to incorporate \gls{QoS}-aware routing by differentiating between flow priorities and sensitivity to service quality, and \textcite{baraszQoSAwareFair2006} consider bounded path lengths and fairness constraints in maximum concurrent flow formulations. These models reflect supply--demand interaction through feasibility constraints and admission control, rather than through a continuous demand response function, and thus differ fundamentally from elasticity-based approaches. Methodological contributions include the limited granularity and limited path heuristic approaches of \textcite{yuan1999a} for solving routing problems with two different \gls{QoS} thresholds, the fully polynomial time approximation scheme of \textcite{huang2010a} for routing problems with multiple \gls{QoS} thresholds, and the sequential heuristic approach of \textcite{bialon2016a} for \gls{QoS}-constrained unsplittable \gls{MCF} problem that iteratively assigns paths and re-routes to resolve capacity violations.

A related but distinct stream of research arises in traffic modeling, where demand--supply interactions are represented through equilibrium frameworks. In these models, such as deterministic or stochastic user equilibrium \parencite{sheffiUrbanTransportationNetworks1984}, users independently choose routes or modes based on perceived travel costs, and demand is distributed across the network accordingly. 
The elastic \gls{TAP} is reviewed by \textcite{gartnerOptimalTrafficAssignment1980a}, where congestion effects are modeled through volume-dependent arc costs. \textcite{gartnerOptimalTrafficAssignment1980b} shows that the elastic \gls{TAP} can be reformulated as an inelastic assignment problem on an expanded network, though dedicated algorithms are still required to solve larger instances. 
The interaction between travel demand and network conditions is typically captured through fixed-point or iterative procedures, such as the method of repeated averages or the method of successive averages, see e.g.,~\textcite{nielsenOptimisationTimetablebasedStochastic2006, richSystemConvergenceTransport2015}, which seek consistency between user choices and resulting travel times, or primal-dual approaches, see e.g.,~\textcite{florian1974a}. While such approaches capture endogenous demand through behavioral responses, they rely on decentralized decision-making and typically ignore explicit capacity allocation and demand admission. 

In summary, the literature shows that supply-dependent demand is increasingly recognized across transportation, logistics, and telecommunications, but existing approaches often either ignore capacity constraints, rely on predefined paths, or treat demand as exogenous or fully served. This motivates a \gls{SAMCF} framework that jointly models routing, capacity allocation, and endogenous demand in a centralized, system-optimal setting.

\vspace{15pt}

\section{Problem Definition}
\label{sec:problem}

This section presents the \gls{SAMCF} problem and formalizes the interaction between routing decisions, service quality, and elastic demand. In Section~\ref{sec:nlp}, we first introduce a nonlinear formulation that extends the classical capacitated MCF problem by incorporating service-dependent demand and an outside alternative. We then illustrate the behavior of the model through a small example in Section~\ref{sec:example}, demonstrating how demand response and capacity constraints jointly influence routing decisions. The main notation is summarized in \autoref{tab:notation}.

\subsection{A Nonlinear Programming Model}
\label{sec:nlp}
We consider a directed network $G = (N, A)$, where $N$ denotes the set of nodes and $A$ denotes the set of directed arcs. Each arc $a\in A$ has a capacity $c_a>0$ and a unit cost $t_a\geq0$. 
Let $\delta^{-}(i)$ and $\delta^{+}(i)$ be the set of incoming and outgoing arcs for node $i$, respectively.
We have a set of commodities $D$.
Each commodity, also termed an OD pair, is defined by a source node $s_d\in N$, a sink node $t_d\in N$, and potential demand $w_d$. 
For each commodity $d\in D$ let variables $f_{da}\geq 0$ be the fractional flow on arc $a$. Flow may be split fractionally across different paths.
In the classical \gls{MCF} problem, the demand $w_d$ is fixed and must be fully routed through the network. The objective is to determine a minimum-cost flow that satisfies flow conservation and capacity constraints.  However, in many applications, demand is not fixed but depends on the level of service provided by the network. For example, in transportation systems, travelers may choose whether to use a public transport service based on expected travel time or cost. 
In telecommunications, users may decide to use a network service based on the expected quality of service (e.g., delay, bandwidth).
In this work, we consider a setting in which demand is responsive to the service level but subject to hard capacity constraints. 
For each commodity $d \in D$, only a fraction $q_d \in [0,1]$ of the potential demand $w_d$ can be routed through the network. 
The realized demand $w_d q_d$ depends both on the service experienced and the available capacity.

We introduce a variable $u_d\geq0$ representing the quality experienced by commodity $d$. Specifically, the service level for commodity $d$ depends on the routing decisions and is modeled as a function of arc flows. 
In general, we write 
\begin{equation}
u_d = s_d(\mathbf{f}),
\end{equation}
where $\mathbf{f}$ is the vector of all arc flows and $s_d(\cdot)$ is the service level function for commodity $d$. 

In this work, we measure the service level for each commodity by the average (or expected) cost experienced by the commodity. That is, the service is given by the total routing cost of the commodity divided by the total flow leaving its source node (provided some flow is routed): 
\begin{equation}
s_d(\mathbf{f}) = \frac{\sum_{a \in A} t_a f_{da}}{\sum_{a \in \delta^{+}(s_d)} f_{da}} \quad \text{if } q_d > 0
\end{equation}
Equivalently, this can be written in bilinear form: 
\begin{equation}
    u_d q_d = \sum_{a \in A} t_a f_{da}
\end{equation}
Similar service-level measures have been used in telecommunications \parencite{leblancPlanningModelsWideArea2002} and transport assignment problems \parencite{cancaIntegratedRailwayRapid2019}. 
We denote by $g=(g_d)_{d\in D}$ the collection of demand functions and model user behavior through a non-increasing function 
\[
p_d = g_d(u_d),
\]
where $p_d \in [0,1]$ represents the fraction of demand that is attracted to the network given service level $u_d$. Importantly, $p_d$ should be interpreted as an upper bound on the demand that {\em can be} attracted to the system rather than the demand that {\em must be} served. 
The routed demand is therefore constrained by
\[
0 \le q_d \le p_d,
\]
allowing for the possibility that not all willing users can be accommodated due to capacity limitations.
This differs from classical choice-based or elastic-demand models, where all attracted demand must be routed. 
Instead, we consider a capacity-constrained setting in which the network operator determines how much of the attracted demand is served and how it is routed.
The demand function $g_d(\cdot)$ can be interpreted as the outcome of a choice process between the network $G$ and an outside competing alternative.
In particular, the fraction of users choosing the network is determined by a logit model:
\begin{equation}
p_d = \frac{1}{1 + e^{\alpha - \beta(\hat{u}_d - u_d)}},
\end{equation}
where parameters $\alpha$ and $\beta$ control the sensitivity of demand to differences in service levels. 
If the network provides a lower cost ($u_d < \hat{u}_d$), a larger share of demand is attracted; otherwise, demand shifts to the competing alternative.
We assume that this alternative is always available and that its service level $\hat{u}_d$ is unaffected by the routing decisions within the network $G$. This assumption is most appropriate when the competing alternative does not share capacity-constrained resources with the network (such as a shared road network) or when the shifted demand is small. We therefore assume that the competing alternative operates on an individual graph $\hat{G}=(N,\hat{A})$ where $\hat{A}$ contains exactly one arc from $s_d$ to $t_d$ per commodity $d$ with unlimited capacity and cost $\hat{u}_d$. 
Additional network-operator-side incentives can be incorporated directly in the objective by adjusting arc costs; for example, per-commodity revenue from routing demand through $G$ can be represented as a negative routing cost (or, equivalently, by increasing the cost of arcs in $\hat{A}$). The value of $\hat{u}_d$ can be kept unchanged in the choice model, so the demand response remains valid.
 
We assume that there is a centralized decision maker that jointly determines, 
i) how much demand is served for each commodity subject to a decentralized choice model, 
ii) how the demand is routed through the network, and 
iii) how the network capacity is allocated across commodities. 
The objective is system-optimal, meaning that the decision-maker minimizes the total routing costs for all demand, whether the commodities are routed through $G$ or through the alternative $\hat{G}$. Given these assumptions, the \gls{SAMCF} problem can be formulated as the following \gls{NLP}:

\begin{align}
\min \quad & \sum_{d \in D} w_d \sum_{a \in A \cup \hat{A}} t_a f_{da} \label{NL:OBJ}\\
\text{s.t.} \quad 
& \sum_{d \in D} w_d f_{da} \leq c_a \qquad && \forall a \in A \label{NL:CAP}\\
& \sum_{a \in \delta^{-}(i)} f_{da} - \sum_{a \in \delta^{+}(i)} f_{da} = b_i^{d} \qquad && \forall i \in N,\ d \in D \label{NL:CONS}\\
& q_{d} = \sum_{a \in \delta^{+}(s_d) \setminus \hat{A}}f_{da}  \qquad && \forall d \in D \label{NL:PTFLOW}\\
& u_d q_d = \sum_{a \in A} t_a f_{da} \qquad && \forall d \in D \label{NL:AVGUTIL}\\
& p_d \leq \frac{1}{1 + e^{\alpha - \beta(\hat{u}_d - u_d)}} \qquad && \forall d \in D \label{NL:PTMODEBOUND} \\
& 0 \leq q_{d}\leq p_d \qquad && \forall d \in D \label{NL:PTCAPBOUND} \\
& 0 \leq f_{da}\leq 1 \qquad && \forall d \in D,\ a \in A \label{NL:FDOM}\\
& 0 \leq p_{d}\leq 1  \qquad && \forall d \in D \label{NL:PDOM} \\
&  0 \leq u_d  \qquad && \forall d \in D \label{NL:UDOM}
\end{align}
The objective function \eqref{NL:OBJ} minimizes total routing costs. Demand routed through the network incurs costs on the utilized arcs,
while unmet demand is assumed to be served by the competing alternative at a
fixed cost $\hat{u}_d$ per unit. 
Constraints \eqref{NL:CAP} enforce arc capacity limits, and constraints~\eqref{NL:CONS} ensure flow conservation for each OD pair, where $b_n^d = -1$ at the origin node $s_d$, $b_n^d = 1$ at the destination node $t_d$, and $b_n^d = 0$ otherwise, for each node $n \in N$. Constraints~\eqref{NL:PTFLOW} define the proportion of demand routed through the network for each commodity.
Constraints~\eqref{NL:AVGUTIL} connect the routing decisions with the experienced service level. Constraints~\eqref{NL:PTMODEBOUND} link service quality to demand through the choice model, and Constraints \eqref{NL:PTCAPBOUND} ensure that routed demand cannot exceed attracted demand. Finally, (\ref{NL:FDOM})--(\ref{NL:UDOM}) define the variable domains.  

The model includes a coupling between flow, service level, capacity, and demand. 
In particular, the service level depends on the flow on the arcs, while demand depends on the service level. 
This results in a nonlinear model characterized by 
i) a bilinear relationship between flow and service level, and ii) a nonlinear relationship between service level and demand. 

We emphasize that, although demand behavior is represented by a choice model, the resulting solution {\em does not} correspond to a classical user equilibrium. 
In particular, the amount of demand attracted to the network need not equal the amount of demand that is served. 
Rather, the model represents a centralized, system-optimal allocation problem in which routing and service levels jointly determine how much demand can be feasibly accommodated.
This interpretation is particularly relevant in strategic planning contexts, where the objective is to assess how much demand a fixed system can capture and how capacity should be allocated across commodities, rather than to model decentralized user behavior exactly.

\begin{table}[htbp!] 
		\centering
        \caption{Notation and terminology}
            \renewcommand{\arraystretch}{0.97}
		\begin{tabular}{ll}
			\toprule
         	\textbf{Sets}& Description\\
            \midrule
            $A$ & Set of arcs\\
            $\hat{A}$ & Set of arcs representing the competing mode \\
            $N$ & Set of nodes\\
            $D$ & Set of commodities\\
            $B$ & Set of PWL breakpoints\\
            \midrule
           \textbf{Parameters} & Description\\
           \midrule
            $c_a$ & Capacity of arc $a$\\
            $t_a$ & Unit cost of arc $a$\\
            $w_d$ & Demand of commodity $d$\\
            $\hat{u}_d$ & Fixed service level for competing mode of commodity $d$\\
            $u_{db}$ & Service level for commodity $d$ at breakpoint $b\in B$ \\
            $p_{db}$ & Fraction of demand captured at breakpoint $b\in B$ for commodity $d$\\
           \midrule
           \textbf{Functions} & Description\\
           \midrule
            $s_d(\cdot)$ & Service level function for commodity $d$\\
            $g_d(\cdot)$ & Demand function for commodity $d$\\
             \midrule
           \textbf{Variables}& Description\\
           \midrule
            $f_{da}$ & Fractional flow of commodity $d$ on arc $a$\\
            $u_{d} $ & Experienced service level for commodity $d$\\
            $p_{d} $ & Fraction of demand attracted for commodity $d$ according to $g_d(u_d)$\\
            $q_{d} $ & Fraction of demand routed for commodity $d$ ($q_d \leq p_d$)\\
            $\lambda_{db}$ & Interpolation weight for breakpoint $b \in B$ of commodity $d$\\   
            $\theta_{db}$ & Auxiliary variables for the linearization of bilinear terms $\lambda_{db}q_d$\\
            \bottomrule
	\end{tabular}
    \label{tab:notation}
\end{table} 

\subsection{An Illustrative Example}
\label{sec:example}
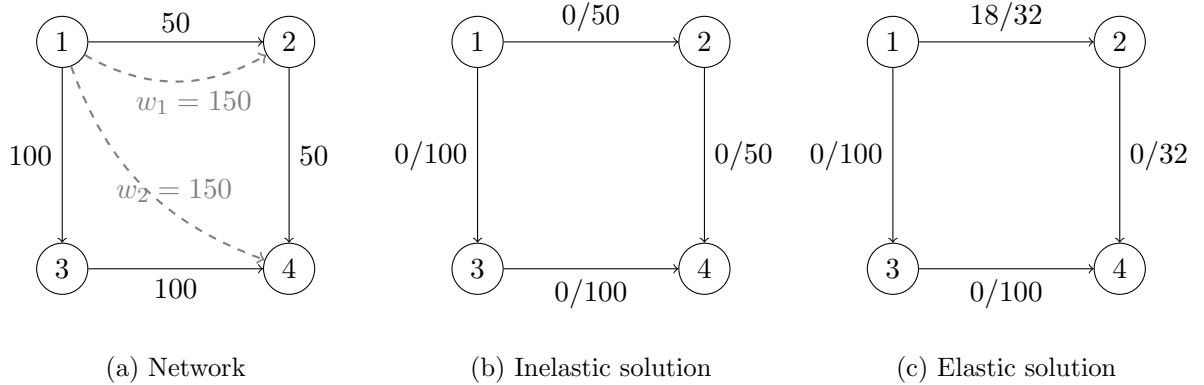
\begin{figure}[htbp]
    \centering
        \centering
          \begin{subfigure}[b]{0.31\textwidth}
          \centering
        \begin{tikzpicture}[baseline=(current bounding box.center)]
            \useasboundingbox (-0.5,0.5) rectangle (3.5,-3.5);
            \node[circle, draw] (1) at (0,0) {1};
            \node[circle, draw] (2) at (3,0) {2};
            \node[circle, draw] (3) at (0,-3) {3};
            \node[circle, draw] (4) at (3,-3) {4};
            \draw[->] (1) -- node[above] {\textcolor{black}{50}} (2);
            \draw[->] (1) -- node[left] {\textcolor{black}{100}} (3);
            \draw[->] (2) -- node[right] {\textcolor{black}{50}} (4);
            \draw[->] (3) -- node[below] {\textcolor{black}{100}} (4);
            \draw[->, dashed, gray, bend right=30, thick] (1) to node[below, pos=0.6] {$w_1 = 150$} (2);
            \draw[->, dashed, gray, bend right=25, thick] (1) to node[above, pos=0.65] {$w_2 = 150$} (4);
        \end{tikzpicture}
        \\[2.5ex]
        \caption{Network}
        \label{fig:tiny-a}
        \end{subfigure}
    \hfill
          \begin{subfigure}[b]{0.31\textwidth}
          \centering
        \begin{tikzpicture}[baseline=(current bounding box.center)]
            \useasboundingbox (-0.5,0.5) rectangle (3.5,-3.5);
            \node[circle, draw] (1) at (0,0) {1};
            \node[circle, draw] (2) at (3,0) {2};
            \node[circle, draw] (3) at (0,-3) {3};
            \node[circle, draw] (4) at (3,-3) {4};
            \draw[->] (1) -- node[above] {\textcolor{black}{0/50}} (2);
            \draw[->] (1) -- node[left] {\textcolor{black}{0/100}} (3);
            \draw[->] (2) -- node[right] {\textcolor{black}{0/50}} (4);
            \draw[->] (3) -- node[below] {\textcolor{black}{0/100}} (4);
            \draw[->, dashed, opacity=0, bend right=30, thick] (1) to node[below, pos=0.6] {} (2);
            \draw[->, dashed, opacity=0, bend right=25, thick] (1) to node[above, pos=0.65] {} (4);
        \end{tikzpicture}
        \\[2.5ex]
        \caption{Inelastic solution}
        \label{fig:tiny-b}
        \end{subfigure}
        \hfill
          \begin{subfigure}[b]{0.31\textwidth}
          \centering        \centering
        \begin{tikzpicture}[baseline=(current bounding box.center)]
            \useasboundingbox (-0.5,0.5) rectangle (3.5,-3.5);
            \node[circle, draw] (1) at (0,0) {1};
            \node[circle, draw] (2) at (3,0) {2};
            \node[circle, draw] (3) at (0,-3) {3};
            \node[circle, draw] (4) at (3,-3) {4};
            \draw[->] (1) -- node[above] {\textcolor{black}{18/32}} (2);
            \draw[->] (1) -- node[left] {\textcolor{black}{0/100}} (3);
            \draw[->] (2) -- node[right] {\textcolor{black}{0/32}} (4);
            \draw[->] (3) -- node[below] {\textcolor{black}{0/100}} (4);
            \draw[->, dashed, opacity=0, bend right=30, thick] (1) to node[below, pos=0.6] {} (2);
            \draw[->, dashed, opacity=0, bend right=25, thick] (1) to node[above, pos=0.65] {} (4);
        \end{tikzpicture}
        \\[2.5ex]
        \caption{Elastic solution}
        \label{fig:tiny-c}
\end{subfigure}
\caption{An illustrative example. \autoref{fig:tiny-a} shows the network where arc labels indicate capacities. Two commodities with demand $w_1=150$ and $w_2=150$ are considered. \autoref{fig:tiny-b} and \autoref{fig:tiny-c} show the optimal solutions of the inelastic and elastic problem, respectively. In the solution figures, each arc label is of the form $f_1/f_2$, where the first value denotes the flow of commodity~1 and the second value denotes the flow of commodity~2.}
    \label{fig:three_networks}
\end{figure}

To demonstrate the effect of considering demand elasticity and a competing alternative, we consider the four-node network in \autoref{fig:tiny-a}. The capacities are shown on each arc, and there are two commodities: commodity 1 from node 1 to 2 with potential demand $w_1=150$, and commodity 2 from node 1 to 4 with potential demand $w_2=150$. The unit cost of each arc is 1. In addition, we assume that the competing alternative has unit costs of $\hat{u}_1 = 3$ and $\hat{u}_2 = 6$. 

\autoref{fig:tiny-b} shows the optimal solution to the inelastic problem with the flow of commodities 1 and 2 shown on the arcs. Since using the competing mode is more costly for commodity 2, the optimal choice is to allocate capacity to commodity 2, while shifting commodity 1 entirely to the competing mode. The resulting objective value is $750$. 
\autoref{fig:tiny-c} shows the solution to the elastic problem with the flow similarly shown on the arcs. We assume a logit function with $\alpha = 0$ and $\beta = 0.5$. Based on the resulting service levels, the fraction of demand attracted is $g_1(1)=0.73$ and $g_2(2)=0.88$, corresponding to at most 110 and 132 units for commodities 1 and 2, respectively. However, due to the capacity constraints, not all the attracted demand can be routed. The optimal solution therefore allocates capacity primarily to commodity 2, routing all of its attracted demand ($q_2 = p_2)$, while only a small fraction of commodity 1 is served ($q_1 \leq p_1)$. The resulting objective value is $786$. 

In this case, the problem is strongly capacitated. This means that while the total amount of demand routed remains the same, the distribution of demand across commodities is different. In other cases, the capacity might be sufficient, and it is the service offered that determines the routing. If capacity constraints were absent, then the inelastic solution would overestimate the demand that could be attracted. 

\section{Solution Methodologies}
\label{sec:methods}

This section presents solution methods for the \gls{SAMCF} problem. We first propose a column generation-based fixed-point heuristic that iteratively couples routing decisions and demand updates. We then introduce two approximation models of \gls{NLP}~\eqref{NL:OBJ}--\eqref{NL:UDOM}: a bilinear piecewise-linear approximation and a fully linearized formulation based on McCormick envelopes. We conclude this section with a discussion on optimality guarantees and the bounds provided by the different approaches. 

\subsection{A Column Generation-Based Fixed-Point Heuristic}
To efficiently solve the \gls{SAMCF}, we propose a heuristic based on an iterative interaction between routing decisions and demand adjustment. 
Each iteration alternates between solving an inelastic assignment for fixed demand and updating demand based on the resulting service levels. 
The approach can be interpreted as a fixed‑point iteration between the supply side (network and capacity) and the demand side (user response), and is repeated until consistency between routing and demand is achieved. 

Our primary goal lies in estimating how much demand is attracted to the network and how it is routed. 
For this, we introduce an iteration-specific potential demand constant $w^{i}_{d} \leq w_d$ per commodity $d \in D$. The algorithm is initialized with an optimistic demand estimate, setting $w^0_d = w_d$, corresponding to the assumption that all potential demand can be routed. 
Routing and demand updates are then performed iteratively until the demand routed through the network $G$ is consistent with the demand implied by the service level. 

\subsubsection{Solving the Inelastic MCF}
For fixed demand values $q_d = w_d^i$, the \gls{SAMCF} problem reduces to a capacitated \gls{MCF} problem in which all demand must be routed. 
Solving this problem directly using an arc-flow formulation quickly becomes intractable for large networks (see e.g.,~\cite{ford1958a}). 
We therefore apply column generation to a path-based reformulation of the arc-flow model.

Let $P_d$ denote the set of all simple $s_d$-$t_d$ paths for commodity $d$, and let $P_d(a)\subseteq P_d$ be the paths that use arc $a \in A \cup \hat{A}$ and $a\in p$ the arcs in path $p$. The cost of path $p$ is $c_p = \sum_{a \in p} t_a$. For each path $p \in P_d$, let $f_{dp} \geq 0$ be the fraction of demand $d$ routed along $p$. 
The path-based master problem is then:
\begin{align}
    \min \quad & \sum_{d \in D} w_d^i \sum_{p \in P_d} c_p \, f_{dp} \label{MP:OBJ}\\
    \text{s.t.} \quad
    & \sum_{d \in D} w_d^i \sum_{p \in P_d(a)}  f_{dp} \leq c_a \qquad && \forall a \in A \label{MP:CAP}\\
    & \sum_{p \in P_d} f_{dp} = 1 \qquad && \forall d \in D \label{MP:DEM}\\
    & f_{dp} \geq 0 \qquad && \forall d \in D,\ p \in P_d \label{MP:NN}
\end{align}
Let $\mu_a \leq 0$ and $\pi_d$ denote the dual variables associated with constraints \eqref{MP:CAP} and \eqref{MP:DEM}, respectively.
Since $|P_d|$ is exponential in the network size, we solve a restricted master problem (RMP) over a subset $\bar{P}_d \subseteq P_d$, initialized with the path using arc $a \in \hat{A}$ for each commodity, i.e., the direct $s_d$-$t_d$ path through the competing mode. 
After solving the RMP, the current dual solution $(\boldsymbol{\mu}, \boldsymbol{\pi})$ is used to identify improving columns via a subproblem.
For each commodity $d$, the reduced cost of a path $p \in P_d$ is:
\begin{equation}
    \bar{c}_{dp} = w_d^i \sum_{a \in p} \bigl(t_a - \mu_a\bigr) - \pi_d. \label{MP:RC}
\end{equation}
A path improves the current solution if and only if $\bar{c}_{dp} < 0$. The subproblem for commodity $d$ finds a minimum reduced-cost path:
\begin{equation}
    p^* \in \arg\min_{p \in P_d}\; \sum_{a \in p} \bigl(t_a - \mu_a\bigr), \label{MP:PP}
\end{equation}
which is a shortest path problem on $G$ with modified arc costs $t_a - \mu_a$. 
If $\bar{c}_{dp^*} < 0$, path $p^*$ is added to the RMP and the procedure repeats. 
The column-generation procedure terminates when no commodity admits a path with negative reduced cost, at which point the RMP solution is optimal for the master problem~\eqref{MP:OBJ}--\eqref{MP:NN}.
The resulting routing $\mathbf{f}$ and adjusted objective
\begin{equation}
    z^i_{\text{MCF}} = \sum_{d \in D} \left(  w_d^i \sum_{p \in \bar{P}_d} c_p f_{dp} + (w_d - w_d^i)\hat{u}_d \right) 
\end{equation}
are recorded, where $(w_d - w_d^i)\hat{u}_d$ accounts for the routing cost of the excluded demand.

This approach allows the inelastic subproblem to be solved efficiently for large networks with many possible paths.
Since the fixed-point heuristic solves the \gls{MCF} problem repeatedly across iterations, the column sets $\bar{P}_d$ accumulated in iteration $i$ are retained and passed as the initial column set for iteration $i+1$. 
Because consecutive iterations differ only in the demand weights $w_d^i$, all of the previously generated paths remain relevant.

\subsubsection{Updating Demand Estimates}
Once the inelastic subproblem is solved, the resulting routing $\mathbf{f}$ is used to compute the experienced service level $u_d = s_d(\mathbf{f})$ for each commodity, and the fraction of demand attracted to the network is updated as $p_d = g_d(u_d)$. We reuse $q_d = \sum_{p \in \bar{P}_d : p \cap \hat{A} = \emptyset} f_{dp}$ to denote the demand routed in the system, here without the restriction that $q_d \leq p_d$. 

The demand estimate for the next iteration is then updated according to the rule
\begin{equation*}
        w^{i+1}_d = \begin{cases}
            \kappa w^i_d, & \ \text{if } q_d = 0, \\
            \kappa q_d w^i_d + (1 - \kappa) p_d w_d, & \text{otherwise}
        \end{cases}
\end{equation*}
where $\kappa \in [0,1]$ is a relaxation parameter controlling how aggressively we update. If no flow is routed for a commodity, the current demand estimate is reduced, reflecting that the network does not provide a sufficiently attractive service for that demand level. Otherwise, two demand estimates are available: the amount of demand that is routed when ignoring the routing quality, and the amount of demand that would be attracted given the experienced service quality. The update rule combines these two estimates using the relaxation parameter $\kappa$, while anchoring the update to the original potential demand $w_d$ to maintain feasibility with respect to maximum demand. This update scheme is conceptually similar to the Method of Successive Averages used in traffic assignment, which combines current estimates with those from previous iterations to stabilize the fixed-point process. Unlike classical MSA approaches based on decentralized user responses, however, the averaging here is applied to endogenous demand estimates within a centralized, capacity-constrained routing framework.

Finally, the routing at any iteration $i$ also admits a feasible solution, and thus a valid upper bound, to the \gls{NLP}~\eqref{NL:OBJ}--\eqref{NL:UDOM}.
For each commodity, we set the repaired demand to $\tilde{w}_d = \min(q_d w_d^i,\, p_d w_d)$, which is the largest amount that is both actually routed through $G$ and consistent with the attracted demand $p_d w_d$.
The repaired objective value is then
\begin{equation}
    z^i_{\text{Repaired}} = \sum_{d \in D} \tilde{w}_d\, u_d + \sum_{d \in D}(w_d - \tilde{w}_d)\hat{u}_d,
    \label{eq:zRepaired}
\end{equation}
and is used in the objective-based termination criterion.

\subsubsection{Termination Criteria}
The heuristic alternates between solving the inelastic \gls{MCF} problem and updating demand estimates until convergence to a fixed point where routing decisions and demand levels are consistent. 
Convergence is assessed using both a demand-based and an objective-based criterion. The algorithm terminates when the change in demand falls below a threshold $\epsilon_w$, when the relative change in objective value is below $\epsilon_{OBJ}$, or when a maximum number of iterations $I$ is reached. 

\subsubsection{Overall Heuristic Framework}
\autoref{main_framework} summarizes the proposed fixed-point heuristic. The output of the heuristic is a routing that is feasible with respect to network capacities and approximately consistent with the demand-response. All algorithm-specific notation is collected in \autoref{tab:notation_algorithm}.

\begin{table}[htbp!] 
		\centering
        \caption{Algorithm-specific notation and terminology}
            \renewcommand{\arraystretch}{0.97}
		\begin{tabular}{ll}
			\toprule
           \textbf{Parameters} & Description\\
           \midrule
            $I$ & Maximum number of iterations\\
            $\epsilon_{OBJ}$ & Objective threshold \\
            $\epsilon_{w}$ & Demand threshold \\
            $\kappa$ & Relaxation factor\\   
            $w^i_d$ & Total demand for commodity $d$ at iteration $i$ \\
            \bottomrule
	\end{tabular}
    \label{tab:notation_algorithm}
\end{table} 

\begin{algorithm}[H]
    \caption{Fixed-point heuristic for the \gls{SAMCF} problem}
    \label{main_framework}
    \KwIn{Graph $G = (N,A \cup \hat{A})$, commodities $D$, a maximum number of iterations $I$, demand and objective value tolerances, $\epsilon_{w}$ and $\epsilon_{OBJ}$, and relaxation factor $\kappa\in[0,1]$}
    \KwResult{Estimate of demand captured (and a routing)} 
    $i \gets 0, \Delta_{w} \gets \infty, \Delta_{OBJ} \gets \infty$\;
    $w^0_d \gets w_d$  \Comment*[r]{Set initial demand estimate}
    \While{$i < I$ \textbf{and} $\Delta_{w} > \epsilon_w$ \textbf{and} $\Delta_{OBJ} > \epsilon_{OBJ}$} {
        $\mathbf{f}, z^i_{\text{MCF}} \gets \textsc{solveMCF}(G, w^{i}_d)$ \Comment*[r]{Solve inelastic problem}
        \ForEach{$d \in D$} {
            $u_d \gets s(\mathbf{f})$\;
            $p_d \gets g(u_d)$\;
            $q_d \gets \sum_{a \in \delta^{+}(s_d)\setminus \hat{A}} f_{da}$\;
            \If{$q_d = 0$}{
                $w^{i+1}_d  \gets \kappa w^{i}_d$ \Comment*[r]{Decrease demand estimate}
                }
            \Else{ 
                $w^{i+1}_d  \gets \kappa q_d w_d^i + (1-\kappa) p_d w_d$ \Comment*[r]{Update demand estimate}
            }
        }
        $i \gets i + 1$\;
        $\Delta_w \gets \sum_{d \in D} |w^{i+1}_d - w^{i}_d|$\; 
        $\Delta_{OBJ} \gets (z^i_{\text{MCF}} - z^i_{\text{Repaired}})/z^i_{\text{Repaired}}$; 
    }
    \Return $\mathbf{f}$\;
    \end{algorithm}

\subsection{A Linear Approximation of the NLP}
We consider two linear approximations of the \gls{NLP}. 
The first replaces the nonlinear demand function in \eqref{NL:PTMODEBOUND} with a piecewise‑linear (PWL) approximation, while retaining the bilinear service–demand coupling in \eqref{NL:AVGUTIL}. The second additionally linearizes this coupling using McCormick envelopes, yielding a fully linear formulation.

\subsubsection{Piecewise‑Linear Approximation of the Demand Function}
For each commodity $d$, the demand function $p_d = g_d(u_d)$ is approximated by a PWL function defined over a finite set of breakpoints $B$. Each breakpoint $b\in B$ is characterized by a service level $u_{db}$ and the corresponding attracted demand fraction $p_{db}=g_d(u_{db})$. The experienced service level $u_d$ and attracted demand $p_d$ are expressed as convex combinations of these breakpoints: 
\begin{align}
    u_d &= \sum_{b \in B} u_{d b} \lambda_{d b} \qquad &&\forall d \in D \\
	p_d &= \sum_{b \in B} p_{d b} \lambda_{d b} \qquad &&\forall d \in D
\end{align}
where $\lambda_{db}$ are interpolation weights satisfying
\begin{align}
   \sum_{b \in B} \lambda_{d b} &= 1 \qquad &&\forall d \in D \\
	\lambda_{db} & \geq 0 \qquad &&\forall d \in D, b \in B.
\end{align}
Additionally, to ensure a valid PWL approximation, we impose an SOS2 constraint on the variables $\lambda_{db}$ for all $b \in B$: 
\begin{align}
    \lambda_{db} &\in \text{SOS2}(1,..,|B|) \qquad && \forall d \in D.
\end{align}
This restricts that at most two adjacent breakpoints can take non-zero values, guaranteeing the demand function is interpolated linearly between neighboring breakpoints. 
By replacing Constraints \eqref{NL:PTMODEBOUND} with the above definitions, we obtain a bilinear program which we will refer to as \pwlBLmodel.

\subsubsection{McCormick Linearization of the Service-Level Function}
In this second approximation, we additionally linearize the bilinear constraint linking flow decisions and service level. Substituting the PWL representation of $u_d$ into constraint \eqref{NL:AVGUTIL} yields: 
$$
u_d q_d = \sum_{a \in A}t_a f_{da} = \sum_{b \in B} u_{db} \lambda_{db} q_d 
$$
which contains bilinear products  $\lambda_{db} q_d$. To linearize these terms, we introduce auxiliary variables $\theta_{db} = \lambda_{db} q_d$. These variables are then approximated using McCormick envelopes, which provide a valid convex relaxation for bounded variables \parencite{mccormickComputabilityGlobalSolutions1976}. The following constraints characterize the relaxation: 
\begin{align}
    \theta_{d b} &\le q_d & \forall d \in D,\ b \in B, \\
    \theta_{d b} &\le \lambda_{d b} & \forall d \in D,\ b \in B, \\
    \theta_{d b} &\ge q_d + \lambda_{d b} - 1 & \forall d \in D,\ b \in B , \\
    \theta_{d b} &\ge 0 & \forall d \in D,\ b \in B.
\end{align}
The service-level constraints \eqref{NL:AVGUTIL} can then be written in linear form as 
\begin{align}
    \sum_{b \in B} u_{db} \theta_{db} = \sum_{a \in A} t_a f_{da} \qquad && \forall d \in D
\end{align}
which yields a \gls{MILP} model that relaxes the BL model. 
We refer to this model as \pwlmodel. 
We note that \pwlBLmodel is an approximation where the accuracy depends on the number and placement of breakpoints and the shape of $g_d$ over the relevant domain. We can control the accuracy by how we set the breakpoints, which we explore in \autoref{sec:results}. 
The \pwlmodel formulation is a convex relaxation and, although tight over bounded domains, may overestimate or underestimate service quality due to fractional auxiliary variables. 
This is a standard trade‑off for obtaining a fully linear formulation.

We note that both the \pwlBLmodel and \pwlmodel are approximations and therefore do not guarantee feasibility for the original \nlmodel. However, a feasible routing and demand estimate can be recovered by computing implied service levels and demand from the selected paths, and then clamping routed demand to $\min\{q_d, p_d\}$ per OD pair, with excess demand redirected to the alternative mode.

\subsection{Optimality and Bounds}

A trivial upper bound for the \gls{SAMCF} problem can be obtained by routing all demand through the alternative arcs in $\hat{A}$, although this solution is typically of poor quality. The proposed heuristic always produces a feasible solution to \gls{NLP} \eqref{NL:OBJ}-\eqref{NL:UDOM}, and therefore provides a valid upper bound. However, as the method is heuristic in nature, no guarantees can be given regarding its optimality gap. The \pwlBLmodel and \pwlmodel formulations are approximations of \gls{NLP} \eqref{NL:OBJ}-\eqref{NL:UDOM}, where the piecewise-linear representation of the demand function may both under- and overestimate the true logit function. As a result, these approximations do not yield valid lower bounds for \gls{NLP} \eqref{NL:OBJ}-\eqref{NL:UDOM}. Moreover, their solutions may not be feasible for the nonlinear formulation, though a feasible solution may be recoverable by repairing the demand, which yields a valid upper bound. Finally, the commercial solver Gurobi can be applied directly to \gls{NLP} \eqref{NL:OBJ}-\eqref{NL:UDOM},  and this is able to provide lower and upper bounds on objective quality if an optimal solution cannot be found within a given time limit.

\subsubsection{An Outer-Approximation of the Demand Function}
We also propose an outer approximation of the demand function. Conceptually, this replaces the original nonlinear demand relation by a conservative PWL representation that yields a relaxation of the nonlinear model and, consequently, valid lower bounds. The construction choices are discussed in the numerical experiments (Section~\ref{sec:numres}).

\section{Computational Experiments}\label{sec:results}

The goal of our study is to evaluate the proposed column-generation-based heuristic (\heurShort) against benchmarks derived from the nonlinear model (\nlmodel) and PWL approximations (\pwlBLmodel and \pwlmodel) on instances motivated by a public transport (PT) application. 

\subsection{A Public Transport Application}
We consider a passenger assignment problem for a public transport system. 
The underlying public transport network (PTN) represents the physical infrastructure in terms of stops and connections, while a \emph{line plan} specifies the operated services, i.e., the selected lines together with their frequencies and capacities.
Given a line plan, the goal is to determine how passengers are distributed across available travel alternatives. 
Commodities correspond to OD pairs, and routing decisions determine the assignment of passengers to paths induced by the line plan. 
In contrast to classical models, demand is elastic and depends on the generalized travel cost experienced in the network.

To more accurately capture the utility of the PT travel options, the line plan is represented as a change-and-go (CNG) network, which is a directed graph that explicitly represents boarding, alighting, in-vehicle travel, and transfers between lines. As described in \autoref{sec:problem}, the alternative mode is represented as direct $s_d-t_d$ arcs with fixed costs. 
We refer to \textcite{schobelLinePlanningMinimal2006} for an introduction to the CNG and to \textcite{hansenExactAlgorithmPublic2026} for details on the specific construction. 

Each OD pair can be routed along multiple public transport paths in the CNG or diverted to a competing mode. The path cost is given by the sum of arc costs and determines the service level experienced by users. 
This service level, in turn, affects the fraction of demand captured by PT through the demand model.

\subsection{Test Instances}\label{sec:instances}
Our test instances are based on line plans designed as part of a previous study on PT line planning. 
The underlying public transport data sets, including the PTN and OD matrices, come from the open-source program LinTim \parencite{schieweLinTimIntegratedEnvironment2024}. The line plans, frequencies, and corresponding CNG graph constructions are based on solutions from \textcite{hansenExactAlgorithmPublic2026} 
Specifically, we use a subset of solutions from Section 6.2 of \textcite{hansenExactAlgorithmPublic2026}, where the networks are designed so that a target set of OD pairs can be routed in the network, avoiding trivial cases in which some OD pairs are completely unserviceable. 

A comparison of the routing assumptions in \textcite{hansenExactAlgorithmPublic2026} and the \gls{SAMCF} formulation is provided in \autoref{apx:EJORobjectives}. The first uses an all-or-nothing assignment in which all demand is routed if an OD-specific threshold is met, meaning changes in demand are more abrupt as all demand for an OD pair switches as the service level changes. In contrast, this work models a smoother, more gradual demand response.

For each line plan, a corresponding CNG graph is constructed based on the selected lines and their operating frequencies. Demand is derived from the LinTim OD matrices, with the Athens and Lower Saxony instances scaled as described in \textcite{hansenExactAlgorithmPublic2026}. Arc costs and capacities are defined consistently with the CNG construction of \textcite{hansenExactAlgorithmPublic2026}, and an operator incentive corresponding to the per-passenger revenue is encoded in the graph.  
Although the original network design problem was formulated using a target subset of 300 OD pairs, our computational experiments consider the complete OD matrix. \autoref{tab:test_problems} summarizes the resulting test instances, reporting for each data set the problem identifier, PTN size (nodes and arcs), CNG size (nodes and arcs), total number of OD pairs, and aggregate travel demand. 

\begin{table}[htbp]
\centering
\caption{Characteristics of the graphs used in the computational experiments.}
\label{tab:test_problems}

\begin{tabular}{@{}llrrrrrr@{}}
\toprule
Data set &
Problem ID &
\multicolumn{2}{c}{PTN} &
\multicolumn{2}{c}{CNG} &
OD Pairs &
Demand \\
\cmidrule(lr){3-4}
\cmidrule(lr){5-6}
\cmidrule(lr){7-7}
\cmidrule(lr){8-8}
&
&
Nodes &
Arcs &
Nodes &
Arcs &
Total &
Total \\
\midrule
Athens &
\TPCodeFromInstance{Athens} &
51 &
52 &
149 &
622 &
2385 &
2957 \\
Grid &
\TPCodeFromInstance{Grid} &
25 &
40 &
143 &
958 &
567 &
2546 \\
Lowersaxony &
\TPCodeFromInstance{Lowersaxony}&
34 &
35 &
111 &
531 &
395 &
1799 \\
Sioux &
\TPCodeFromInstance{Sioux} &
24 &
38 &
109 &
636 &
552 &
4115 \\
\bottomrule
\end{tabular}
\end{table}

\begin{figure}[h]
\centering
\begin{subfigure}{0.24\textwidth}
\centering
\includegraphics[width=\linewidth]{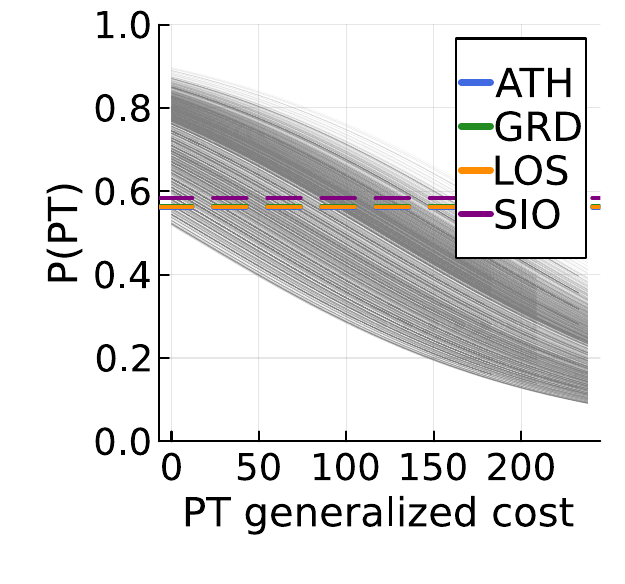}
\caption{$\beta = 0.01$}
\end{subfigure}\hfill
\begin{subfigure}{0.24\textwidth}
\centering
\includegraphics[width=\linewidth]{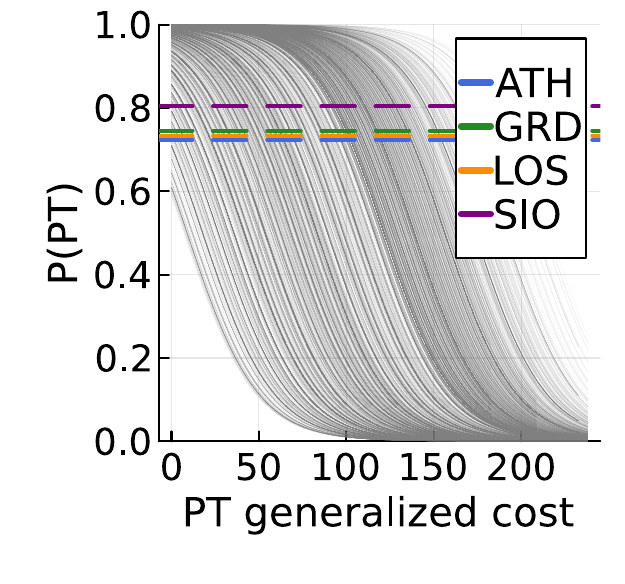}
\caption{$\beta = 0.05$}
\end{subfigure}\hfill
\begin{subfigure}{0.24\textwidth}
\centering
\includegraphics[width=\linewidth]{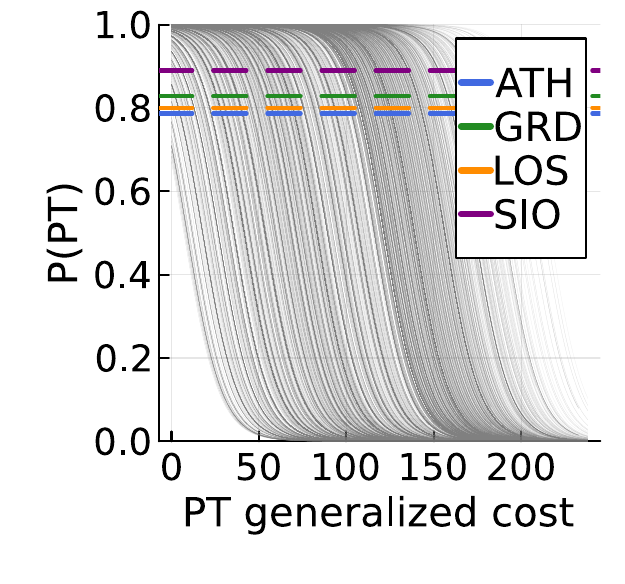}
\caption{$\beta = 0.1$}
\end{subfigure}\hfill
\begin{subfigure}{0.24\textwidth}
\centering
\includegraphics[width=\linewidth]{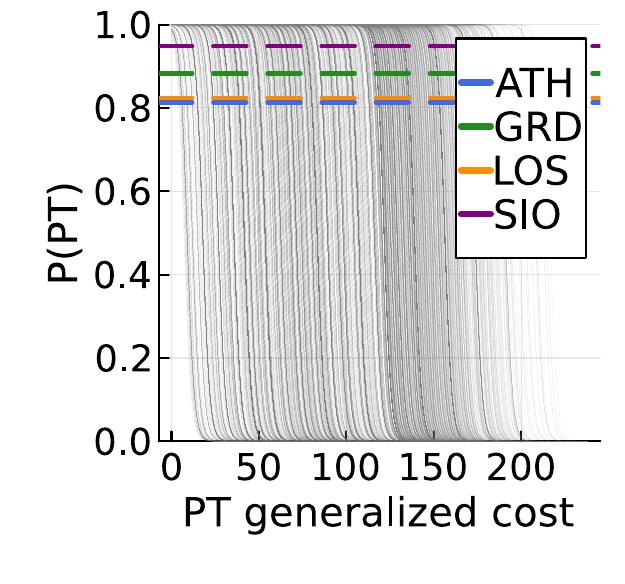}
\caption{$\beta = 0.5$}
\end{subfigure}
\caption{Logit functions for different values of $\beta$ across all OD pairs in the test instances. The dashed lines indicate the corresponding free-flow mode share, which is the fraction of demand that would be attracted to the network if all OD pairs could choose their shortest path without capacity constraints.}
\label{fig:logit_plots}
\end{figure}

We consider a range of logit parameters $\beta$ to capture different levels of sensitivity to service quality. 
The logit parameters $\alpha$ and $\beta$ control the shape of the demand function: $\alpha$ shifts the function horizontally and represents a baseline preference not affected by changes in service quality, while $\beta$ determines the sensitivity to changes in service quality. 
\autoref{fig:logit_plots} visualizes the demand functions \eqref{NL:PTMODEBOUND} for the different $\beta$-settings across all test problems. The figure shows the PT mode share as a function of PT utility, i.e., the service-level defined by \eqref{NL:AVGUTIL}. The average mode share for the ``free-flow'' PT costs (no capacity constraints) is shown for each test problem. 
As seen in \autoref{fig:logit_plots}, larger values of $\beta$ correspond to very sensitive demand functions in which small changes in the PT utility lead to large changes, i.e., almost an all-or-nothing behavior. In contrast, smaller values of $\beta$ correspond to more gradual demand functions, for example, for $\beta = 0.01$ the demand function is almost flat.

For each network and each value of $\beta$, we generate four instances of varying size by including the OD pairs that account for $40\%$, $60\%$, $80\%$, and $100\%$ of total demand, yielding a total of 64 instances across the four networks.

\subsection{Numerical Results}\label{sec:numres}

We evaluate the proposed \heurShort against the nonlinear formulation (\nlmodel) and its piecewise-linear approximations, namely the bilinear model (\pwlBLmodel) and the fully linearized model (\pwlmodel), on the 64 test instances described.
For the \heurShort, the termination criteria are set to a maximum of $I=100$ iterations, $\epsilon_w=0.1$, and a relative objective improvement threshold of $\epsilon_{OBJ}=0.1\%$. The relaxation factor is set to $\kappa=0.5$.
For the \pwlmodel and \pwlBLmodel models, we report recovered objective values and demand estimates, as this allows a consistent comparison across methods. For both \pwlBLmodel and \pwlmodel, the demand functions are approximated using adaptively generated breakpoints with a target sum of squared errors (SSE) of $0.01$. As shown in \autoref{sec:approx}, this setting provides the best trade-off between approximation accuracy and computational effort.
The heuristic is implemented in Julia~1.12.1 and executed using four threads. All optimization models are formulated in JuMP and solved with Gurobi~13.0.0 on a machine with an Intel Xeon Gold 6226R processor and 20 GB of RAM. A time limit of 10 minutes is imposed for each instance. We additionally investigate bound quality using a two-hour time limit for the \nlmodel, \pwlBLmodel, and \pwlmodel, where we use an outer-approximation approach to construct the PWL approximation, ensuring that the resulting formulation provides a valid relaxation and therefore a valid lower bound. The outer approximation is constructed by first finding the optimal breakpoint positions for a target SSE of 0.001, then shifting all points vertically by the largest error between the approximation and $g$.

The analysis first focuses on the comparison between \heurShort, \pwlmodel, and \nlmodel, as the study in \autoref{sec:approx} demonstrated that \pwlmodel consistently provides tighter solutions and better computational performance than \pwlBLmodel. 

\begin{table}[htbp]
\caption{Instance-level comparison of the heuristic (\heurShort), \pwlmodel, and \nlmodel for $\beta = 0.05$.}\label{tab:beta0p05}
\resizebox{\textwidth}{!}{%
\begin{tabular}{@{}r 
            @{\hspace{8mm}} rrr
            @{\hspace{8mm}}
            rrr
            @{\hspace{8mm}}
            rr
            @{\hspace{8mm}}
            rrr@{}}
\toprule
  & \multicolumn{3}{c}{Comp. time (s)} & \multicolumn{3}{c}{Objective value} & \multicolumn{2}{c}{Lower bound} & \multicolumn{3}{c}{PT demand} \\ \cmidrule(lr){2-4} \cmidrule(lr){5-7} \cmidrule(lr){8-9} \cmidrule(lr){10-12}   
  \#ODs & \heurShort & \nlmodel & \pwlmodel & \heurShort & \nlmodel & \pwlmodel & \nlmodel & \pwlmodel & \heurShort & \nlmodel & \pwlmodel \\ 
 \addlinespace[3pt]
 \multicolumn{12}{l}{\textbf{\TPCodeFromInstance{Athens}}} \\
172 & 0.1 & (13.01\%) & 71.75 & 63,064.39 & 63,688.20 & 63,069.62 & 55,399.71 & 60,904.76 & 814.71 & 808.02 & 814.64 \\
453 & 0.23 & -- & (0.05\%) & 104,676.39 & -- & 104,754.57 & 93,038.01 & 101,382.24 & 1,211.08 & -- & 1,204.65 \\
982 & 0.52 & -- & (0.11\%) & 149,979.92 & -- & 150,093.91 & 134,357.10 & 145,378.08 & 1,611.62 & -- & 1,600.96 \\
2,385 & 1.1 & -- & -- & 203,752.71 & -- & -- & 184,411.72 & -- & 2,019.69 & -- & -- \\
\addlinespace[3pt]
 \multicolumn{12}{l}{\textbf{\TPCodeFromInstance{Grid}}} \\
39 & 0.05 & (7.58\%) & 9.01 & 71,071.39 & 71,071.37 & 71,072.59 & 65,685.24 & 68,870.87 & 768.71 & 768.72 & 768.70 \\
101 & 0.08 & (9.61\%) & 127.27 & 104,027.46 & 104,027.42 & 104,028.70 & 94,025.25 & 100,839.13 & 1,136.70 & 1,136.70 & 1,136.68 \\
228 & 0.1 & (29.69\%) & 407.87 & 146,586.59 & 187,181.09 & 146,589.60 & 131,606.52 & 142,129.10 & 1,518.76 & 1,115.55 & 1,518.72 \\
567 & 0.29 & -- & -- & 201,139.20 & -- & -- & 182,293.82 & -- & 1,898.09 & -- & -- \\
\addlinespace[3pt]
 \multicolumn{12}{l}{\textbf{\TPCodeFromInstance{Lowersaxony}}} \\
12 & 0.04 & 9.26 & 0.36 & 26,471.49 & 26,471.47 & 26,471.48 & 26,471.47 & 25,097.36 & 423.60 & 423.60 & 423.60 \\
31 & 0.05 & (3.92\%) & 1.85 & 42,006.21 & 42,006.19 & 42,009.80 & 40,359.87 & 39,886.00 & 620.86 & 620.87 & 620.82 \\
78 & 0.07 & (9.39\%) & 5.99 & 60,511.29 & 60,511.27 & 60,518.59 & 54,827.90 & 57,706.14 & 852.99 & 852.99 & 852.90 \\
395 & 0.15 & (26.04\%) & (0.45\%) & 83,929.54 & 101,140.34 & 83,943.09 & 74,799.07 & 80,035.37 & 1,106.17 & 892.41 & 1,106.00 \\
\addlinespace[3pt]
 \multicolumn{12}{l}{\textbf{\TPCodeFromInstance{Sioux}}} \\
61 & 0.05 & (11.14\%) & 7.1 & 82,887.36 & 82,887.33 & 82,894.79 & 73,650.92 & 80,273.57 & 1,240.67 & 1,240.67 & 1,240.58 \\
116 & 0.07 & (11.32\%) & 41.25 & 139,181.67 & 139,181.63 & 139,207.28 & 123,419.60 & 135,211.64 & 1,893.54 & 1,893.54 & 1,893.26 \\
209 & 0.1 & (11.01\%) & 443.7 & 199,087.56 & 199,087.49 & 199,123.72 & 177,168.11 & 193,708.20 & 2,555.71 & 2,555.71 & 2,555.31 \\
552 & 0.22 & -- & (0.09\%) & 268,568.11 & -- & 268,617.94 & 240,345.60 & 261,287.81 & 3,214.96 & -- & 3,214.41 \\
\bottomrule
\end{tabular}
}
\end{table}

\autoref{tab:beta0p05} reports detailed results for $\beta=0.05$. Instance-level results are reported in \autoref{apx:detailed_tables} for the other values of $\beta$. For each instance, we report the solution time, the objective value, and the captured PT demand (number of passengers routed in the network). For the optimization models, we additionally report the lower bound at termination. Entries shown in parentheses are the optimality gaps reported by Gurobi for instances that time out and ``--`` indicates that no feasible solution was found within the time limit or the program ran out of memory. The gap reported for the \pwlmodel models refers to the optimality gap reported by Gurobi and not the gap between the lower bound and the recovered objective value.

The \heurShort solves all 16 instances in at most $1.1$ seconds, while the optimization models reach the time limit on the largest instances of every network. The \nlmodel times out on 15 of 16 instances and reports gaps between $4$--$30\%$.
However, the objective reported by the \nlmodel closely matches that of the \heurShort in 8 instances, indicating that it has found high-quality solutions but struggles to prove optimality. 

The \pwlmodel scales better, solving 10 of 16 instances to optimality and finding solutions within $0.45\%$ of the reported lower bound on four additional instances. The full $100\%$-demand instances of \TPCodeFromInstance{Athens} and \TPCodeFromInstance{Grid} could not be solved within the time limit, as they exceeded the available memory limit of 20 GB.  

When the models do solve within reasonable bounds, their objective values and captured demand closely match those of the \heurShort. Across the 14 instances solved by the \pwlmodel within the time limit, the \heurShort obtains solutions that are within $0.08\%$ of the model objective value, while the difference in routed demand remains below $0.67\%$.
This strongly suggests that the \heurShort is finding solutions of high quality. By contrast, when the \nlmodel terminates with a larger gap, it tends to substantially underestimate captured demand: 
on \TPCodeFromInstance{Grid} with 228 OD pairs it reports an objective $27.69\%$ above the heuristic and routes $26.55\%$ fewer passengers in the PT network; on \TPCodeFromInstance{Lowersaxony} with 395 OD pairs the corresponding figures are $20.51\%$ and $19.32\%$. 
However, the \heurShort is not always the best-performing method. 
For some instances, the heuristic reports slightly higher objective values than the \nlmodel or \pwlmodel models, but deviations are less than 0.7. 

\begin{table}[htbp]
\centering
\small
\setlength{\tabcolsep}{5pt}
\renewcommand{\arraystretch}{1.08}
\caption{Aggregate performance by demand-sensitivity level $\beta$ and model. Solved indicates the number of instances solved within the time limit (out of 16). Mean reported optimality gap is shown for optimization models. Mean gap to best is computed relative to the best objective found per instance. PT demand statistics are computed on the common subset of instances for which all considered models produced a feasible solution. Mean PT demand reports the average number of passengers captured. Maximum and minimum demand differences indicate, respectively, the worst and closest-to-best relative differences in PT demand compared with the best-performing model for the same instance.}
\begin{tabular}{llrrrrrr}
\toprule
&  & \# Solved & Mean time & Mean reported & Mean gap to & \multicolumn{2}{c}{PT demand} \\
 & & (of 16) & (s) & opt. gap (\%) & best (\%) & max dif.(\%) & min dif.(\%)\\
\midrule
\multirow{4}{*}{$\beta=0.01$} & \heurShort & 16 & 0.23 &  & 0.0 & -0.0001 & 0.0 \\
 & \nlmodel & 10 & 540.55 & 5.55 & 0.12 & -7.11 & 0.0 \\
 & \pwlBLmodel & 11 & 334.12 & 0.52 & 0.45 & -3.41 & -0.09 \\
& \pwlmodel & 15 & 129.84 & 0.04 & 0.02 & -1.13 & 0.0 \\
\addlinespace[6pt]
\multirow{4}{*}{$\beta=0.05$} & \heurShort & 16 & 0.20 &  & 0.0 & -0.0001 & 0.0 \\
& \nlmodel & 11 & 546.30 & 12.07 & 4.47 & -26.55 & 0.0 \\
 & \pwlBLmodel & 14 & 417.22 & 3.22 & 3.06 & -33.12 & -0.36 \\
 & \pwlmodel & 14 & 251.15 & 0.05 & 0.02 & -0.66 & -0.0 \\
\addlinespace[6pt]
\multirow{4}{*}{$\beta=0.1$} & \heurShort & 16 & 0.21 &  & 0.0 & -0.0002 & 0.0 \\
 & \nlmodel & 10 & 542.84 & 8.45 & 1.61 & -7.38 & 0.0 \\
& \pwlBLmodel & 13 & 438.80 & 4.93 & 5.38 & -48.74 & -0.55 \\
& \pwlmodel & 14 & 238.30 & 0.16 & 0.14 & -1.87 & -0.03 \\
\addlinespace[6pt]
\multirow{4}{*}{$\beta=0.5$} & \heurShort & 16 & 0.19 &  & 0.0 & -0.0006 & 0.0 \\
 & \nlmodel & 7 & 465.32 & 4.93 & 3.22 & -32.59 & 0.0 \\
 & \pwlBLmodel & 13 & 504.41 & 4.37 & 6.12 & -87.72 & -0.35 \\
& \pwlmodel & 15 & 132.22 & 0.04 & 0.18 & -0.25 & -0.06 \\
\bottomrule
\end{tabular}
\label{tab:agg_table}
\end{table}

While \autoref{tab:beta0p05} reports instance-level outcomes, \autoref{tab:agg_table} provides an aggregated view across test problems and $\beta$-values. 
For each $\beta$ and method, we report the number of instances solved within the time limit, the mean solution time, the mean reported optimality gap, the mean ``gap to best'' computed against the best objective value found across all methods on each instance, and the maximum and minimum percentage differences from the PT demand captured in the best solution found by any model. As the optimization models could not solve all instances, the mean is computed only on the available values. 

Three observations stand out. 
First, the \heurShort dominates on robustness and runtime: it solves all 64 instances with a mean runtime of $0.21$ seconds, more than three orders of magnitude faster than the models, and it achieves a mean gap to best of zero across every setting of $\beta$. 
The \nlmodel solves only $38$ of $64$ instances within the time limit.
The \pwlmodel is the most competitive benchmark, solving $58$ instances with mean reported gaps below $0.2\%$, but its mean solution time is still around $3$ minutes.
Second, on solution quality, the \pwlmodel remains very close to the best known solutions for all values of $\beta$, with mean gaps to best less than 0.18\%. By contrast, the \pwlBLmodel and \nlmodel are further from the best solution.  
Third, the heuristic consistently achieves the highest average demand across all values of $\beta$, but the \pwlmodel produces nearly identical demand levels: the relative difference in demand to the best solution remains below $0.2\%$ for all settings.
The \nlmodel and \pwlBLmodel models generally attract fewer passengers than the \heurShort, and as $\beta$ increases, the differences become more pronounced. 
For $\beta=0.5$, the heuristic attains the highest average demand (1525 passengers), while the optimization models capture between 1379 and 1522 passengers on average. The worst-case demand losses relative to the best-performing solution on the same instance reach $32.59\%$ for the \nlmodel, $87.72\%$ for the \pwlBLmodel, while the \pwlmodel remains very close to the best-known solution. 
This shows that when the \nlmodel and \pwlBLmodel models fail to close the gap within the time limit, their incumbent solutions tend to under-route the network and divert demand to the alternative mode. 

Taken together, these results indicate that the \heurShort consistently finds solutions of comparable quality to the best optimization-based benchmark while being orders of magnitude faster and substantially more robust across instance sizes and demand assumptions.

\subsubsection{Examining Bound Quality}
\label{sec:twohour}
\begin{figure}[htbp]
    \centering
\begin{tikzpicture}
\begin{axis}[
    ybar,
    bar width=10pt,
    width=12cm,
    height=7cm,
    xlabel={$\beta$},
    ylabel={Mean relative gap (\%)},
    xlabel style={text=black},
    ylabel style={text=black},
    tick label style={text=black},
    symbolic x coords={0.01,0.05,0.1,0.5},
    xtick=data,
    grid=major,
    legend style={at={(0.5,1.04)}, anchor=south, legend columns=3, text=black, draw=black},
]
\addplot+[
    draw=none,
    fill=blue!70,
    bar shift=-10pt,
    error bars/.cd,
    y dir=both,
    y explicit,
    error bar style={black},
    error mark options={black},
] table[
    x=beta,
    y=nl_mean_gap,
    y error plus=nl_err_high,
    y error minus=nl_err_low,
    col sep=comma
] {fig9b_tikz_wide.csv};
\addlegendentry{\nlmodel}
\addplot+[
    only marks,
    mark=none,
    bar shift=-10pt,
    forget plot,
    nodes near coords,
    point meta=explicit symbolic,
    nodes near coords={\pgfplotspointmeta},
    every node near coord/.append style={font=\scriptsize, anchor=south, yshift=2pt, text=black},
] table[
    x=beta,
    y expr=\thisrow{nl_mean_gap}+\thisrow{nl_err_high},
    col sep=comma
] {fig9b_tikz_wide.csv};
\addplot+[
    draw=none,
    fill=green!60!black,
    bar shift=0pt,
    error bars/.cd,
    y dir=both,
    y explicit,
    error bar style={black},
    error mark options={black},
] table[
    x=beta,
    y=bl_mean_gap,
    y error plus=bl_err_high,
    y error minus=bl_err_low,
    col sep=comma
] {fig9b_tikz_wide.csv};
\addlegendentry{\pwlBLmodel}
\addplot+[
    only marks,
    mark=none,
    bar shift=0pt,
    forget plot,
    nodes near coords,
    point meta=explicit symbolic,
    nodes near coords={\pgfplotspointmeta},
    every node near coord/.append style={font=\scriptsize, anchor=south, yshift=2pt, text=black},
] table[
    x=beta,
    y expr=\thisrow{bl_mean_gap}+\thisrow{bl_err_high},
    col sep=comma
] {fig9b_tikz_wide.csv};
\addplot+[
    draw=none,
    fill=orange!85!black,
    bar shift=10pt,
    error bars/.cd,
    y dir=both,
    y explicit,
    error bar style={black},
    error mark options={black},
] table[
    x=beta,
    y=mc_mean_gap,
    y error plus=mc_err_high,
    y error minus=mc_err_low,
    col sep=comma
] {fig9b_tikz_wide.csv};
\addlegendentry{\pwlmodel}
\addplot+[
    only marks,
    mark=none,
    bar shift=10pt,
    forget plot,
    nodes near coords,
    point meta=explicit symbolic,
    nodes near coords={\pgfplotspointmeta},
    every node near coord/.append style={font=\scriptsize, anchor=south, yshift=2pt, text=black},
] table[
    x=beta,
    y expr=\thisrow{mc_mean_gap}+\thisrow{mc_err_high},
    col sep=comma
] {fig9b_tikz_wide.csv};
\end{axis}
\end{tikzpicture}
    \caption{Relative gap between the lower bounds from \nlmodel, \pwlBLmodel, and \pwlmodel and the objective value obtained by \heurShort. The bar shows the average gap across the instances considered, and the error bars indicate the minimum and maximum gap across instances.}
    \label{fig:bounds}
\end{figure}
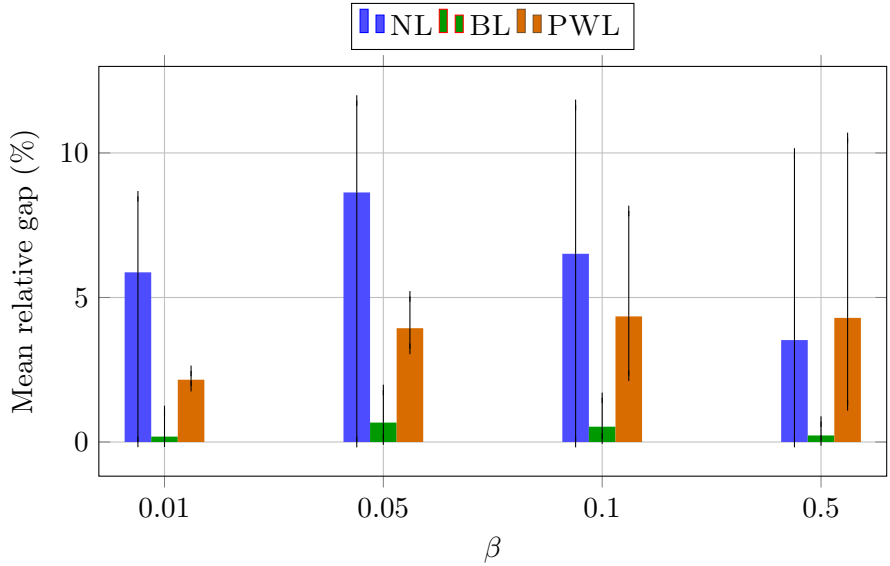

We observe that both the \nlmodel and \pwlBLmodel struggle to close the gap between the incumbent solution and the lower bound within the 10-minute time limit. 
Furthermore, while the recovered solutions from the \pwlmodel are close to the solutions from \heurShort, the lower bounds are relatively weak, with gaps of up to $9.85\%$ to the best feasible solution. 
To complement the previous analysis, we therefore extend the computational experiments by increasing the time limit to two hours for the \nlmodel, \pwlBLmodel, and \pwlmodel.
For the latter two models, the demand approximation is constructed using the outer approximation scheme, such that the resulting formulations provide valid relaxations of the nonlinear demand constraints and, therefore, yield valid lower bounds.
To evaluate the tightness of the bounds around the solutions reported by the \heurShort, we report the relative gap between the lower bounds obtained by the optimization models and the solution found by the \heurShort. 
The relative gap is computed as $\text{Gap}(\%) = \frac{z_{\heurShort} - z_{\text{LB}}}{z_{\heurShort}} \times 100$ where \(z_{\heurShort}\) denotes the objective value of the solution obtained by the \heurShort and \(z_{\text{LB}}\) denotes the lower bound provided by the optimization model.
\autoref{fig:bounds} shows the average gap across instances for each $\beta$-setting, with error bars indicating the minimum and maximum gap observed across instances. The instance-level results underlying the figure are reported in \autoref{apx:twohour}.

The analysis clearly shows that the \pwlBLmodel using the outer-approximation provides the tightest lower bounds overall: its mean gap to the \heurShort solution stays below 1\% for every setting of $\beta$ with relatively low variation across instances. By comparison, the \nlmodel yields weak bounds even with the additional solve time, with an average gap roughly between 3.5\% and 8.5\%. Moreover, increasing the time limit from 10 minutes to 2 hours yields only a marginal improvement in the obtained bounds, reducing the average gap by just 0.5 percentage points. 
The \pwlmodel obtains lower bounds within 2.2\% to 4.3\% of the \heurShort solutions, on average. However, the variability is large, especially as $\beta$ increases. 
This indicates that the McCormick relaxation becomes less accurate when small differences in service level cause large changes in attracted demand.
Overall, this indicates that the bilinear formulation is capable of producing competitive lower bounds when sufficient computational time is available. At the same time, the small average gaps provide further evidence that the \heurShort solutions are close to the best provable bounds.

\subsubsection{Approximation Trade-Offs} \label{sec:approx}
The strength of the \pwlBLmodel and \pwlmodel formulations depends on the number and placement of breakpoints, which jointly determine approximation quality and computational tractability. Higher breakpoint resolution improves the fidelity of $g$, but increases model size and solve time.
To explore this trade-off, we consider two breakpoint strategies: \equidist, using 5, 10, and 20 evenly spaced breakpoints, and \opt, where breakpoints are chosen to achieve a certain approximation error (SSE of 1, 0.1, and 0.01) relative to each demand function in $g$. 

We evaluate recoverability in terms of recovered solution quality, computation time, and robustness across the 64 instances. 
To enable comparisons across instances, recovered objective values are normalized relative to the solution obtained by \heurShort, which found feasible solutions for all instances. A normalized value of 100, therefore, corresponds to the \heurShort baseline, while a value of 101 indicates an average objective that is 1\% higher (worse). 

\begin{figure}[htbp]
\centering
\begin{tikzpicture}
\begin{axis}[
  width=0.85\linewidth,
  height=0.60\linewidth,
  xlabel={Mean comp. time (s)},
  ylabel={Mean normalized obj (\heurShort $\to$ 100)},
  grid=major,
  grid style={line width=0.3pt, draw=gray!30},
  tick label style={font=\small},
  label style={font=\small},
  legend style={at={(0.02,0.98)}, anchor=north west, font=\footnotesize,
    draw=gray!60, fill=white, fill opacity=0.85},
    legend cell align=left,
  clip=false,
  enlarge x limits=0.15,
  enlarge y limits=0.10,
]
\draw[dashed, gray!60, line width=1.6pt]
  (axis cs:50,100) -- (axis cs:600,100);

  \addplot[
    only marks, mark=*, mark size=8.32pt,
    mark options={fill=colBL, fill opacity=0.75, draw=white, line width=0.5pt},
    point meta=explicit symbolic,
    nodes near coords, nodes near coords style={font=\scriptsize, anchor=west, xshift=4pt},
    forget plot,
  ] coordinates {(418.1859, 106.0747) [B=5]};

  \addplot[
    only marks, mark=*, mark size=8.36pt,
    mark options={fill=colBL, fill opacity=0.75, draw=white, line width=0.5pt},
    point meta=explicit symbolic,
    nodes near coords, nodes near coords style={font=\scriptsize, anchor=west, xshift=-35pt},
    forget plot,
  ] coordinates {(500.3012, 105.7991) [B=10]};

  \addplot[
    only marks, mark=*, mark size=8.27pt,
    mark options={fill=colBL, fill opacity=0.75, draw=white, line width=0.5pt},
    point meta=explicit symbolic,
    nodes near coords, nodes near coords style={font=\scriptsize, anchor=west, xshift=4pt},
    forget plot,
  ] coordinates {(531.9291, 105.9136) [B=20]};

  \addplot[
    only marks, mark=square*, mark size=8.55pt,
    mark options={fill=colBL, fill opacity=0.75, draw=white, line width=0.5pt},
    point meta=explicit symbolic,
    nodes near coords, nodes near coords style={font=\scriptsize, anchor=west, xshift=4pt},
    forget plot,
  ] coordinates {(427.3878, 103.8697) [SSE=0.01]};

  \addplot[
    only marks, mark=square*, mark size=8.36pt,
    mark options={fill=colBL, fill opacity=0.75, draw=white, line width=0.5pt},
    point meta=explicit symbolic,
    nodes near coords, nodes near coords style={font=\scriptsize, anchor=west, xshift=4pt},
    forget plot,
  ] coordinates {(346.1952, 100.8418) [SSE=0.1]};

  \addplot[
    only marks, mark=square*, mark size=8.51pt,
    mark options={fill=colBL, fill opacity=0.75, draw=white, line width=0.5pt},
    point meta=explicit symbolic,
    nodes near coords, nodes near coords style={font=\scriptsize, anchor=west, xshift=4pt},
    forget plot,
  ] coordinates {(345.8176, 102.2630) [SSE=1]};

  \addplot[
    only marks, mark=*, mark size=8.96pt,
    mark options={fill=colMcCormick, fill opacity=0.75, draw=white, line width=0.5pt},
    point meta=explicit symbolic,
    nodes near coords, nodes near coords style={font=\scriptsize, anchor=west, xshift=4pt},
    forget plot,
  ] coordinates {(222.6825, 101.8687) [B=5]};

  \addplot[
    only marks, mark=*, mark size=8.6pt,
    mark options={fill=colMcCormick, fill opacity=0.75, draw=white, line width=0.5pt},
    point meta=explicit symbolic,
    nodes near coords, nodes near coords style={font=\scriptsize, anchor=west, xshift=4pt,yshift=3pt},
    forget plot,
  ] coordinates {(325.6529, 100.0595) [B=10]};

  \addplot[
    only marks, mark=*, mark size=8.6pt,
    mark options={fill=colMcCormick, fill opacity=0.75, draw=white, line width=0.5pt},
    point meta=explicit symbolic,
    nodes near coords, nodes near coords style={font=\scriptsize, anchor=west, xshift=-33pt,yshift=-6pt},
    forget plot,
  ] coordinates {(314.7085, 100.0241) [B=20]};

  \addplot[
    only marks, mark=square*, mark size=8.87pt,
    mark options={fill=colMcCormick, fill opacity=0.75, draw=white, line width=0.5pt},
    point meta=explicit symbolic,
    nodes near coords, nodes near coords style={font=\scriptsize, anchor=west, xshift=5pt,yshift=4pt},
    forget plot,
  ] coordinates {(186.1988, 100.0887) [SSE=0.01]};

  \addplot[
    only marks, mark=square*, mark size=9.0pt,
    mark options={fill=colMcCormick, fill opacity=0.75, draw=white, line width=0.5pt},
    point meta=explicit symbolic,
    nodes near coords, nodes near coords style={font=\scriptsize, anchor=west, xshift=-13pt,yshift=12pt},
    forget plot,
  ] coordinates {(143.8114, 100.2044) [SSE=0.1]};

  \addplot[
    only marks, mark=square*, mark size=9.0pt,
    mark options={fill=colMcCormick, fill opacity=0.75, draw=white, line width=0.5pt},
    point meta=explicit symbolic,
    nodes near coords, nodes near coords style={font=\scriptsize, anchor=west, xshift=-13pt,yshift=12pt},
    forget plot,
  ] coordinates {(75.5787, 100.8540) [SSE=1]};

\addlegendimage{only marks, mark=*, mark size=3pt,
  mark options={fill=colBL, draw=white}}
\addlegendentry{\pwlBLmodel (\equidist)}
\addlegendimage{only marks, mark=square*, mark size=3pt,
  mark options={fill=colBL, draw=white}}
\addlegendentry{\pwlBLmodel (\opt)}
\addlegendimage{only marks, mark=*, mark size=3pt,
  mark options={fill=colMcCormick, draw=white}}
\addlegendentry{\pwlmodel (\equidist)}
\addlegendimage{only marks, mark=square*, mark size=3pt,
  mark options={fill=colMcCormick, draw=white}}
\addlegendentry{\pwlmodel (\opt)}
\end{axis}
\end{tikzpicture}
\caption{Trade-off between computation time and solution quality for different breakpoint strategies. The x-axis shows the average computation time in seconds, the y-axis shows the average objective value of the recovered, \nlmodel-feasible solutions obtained from the linear (square markers) and bilinear (round markers) approximations. 
The objective values have been normalized according to \heurShort $= 100$.}
\label{fig:tradeoff_all_qos_baseline}
\end{figure}
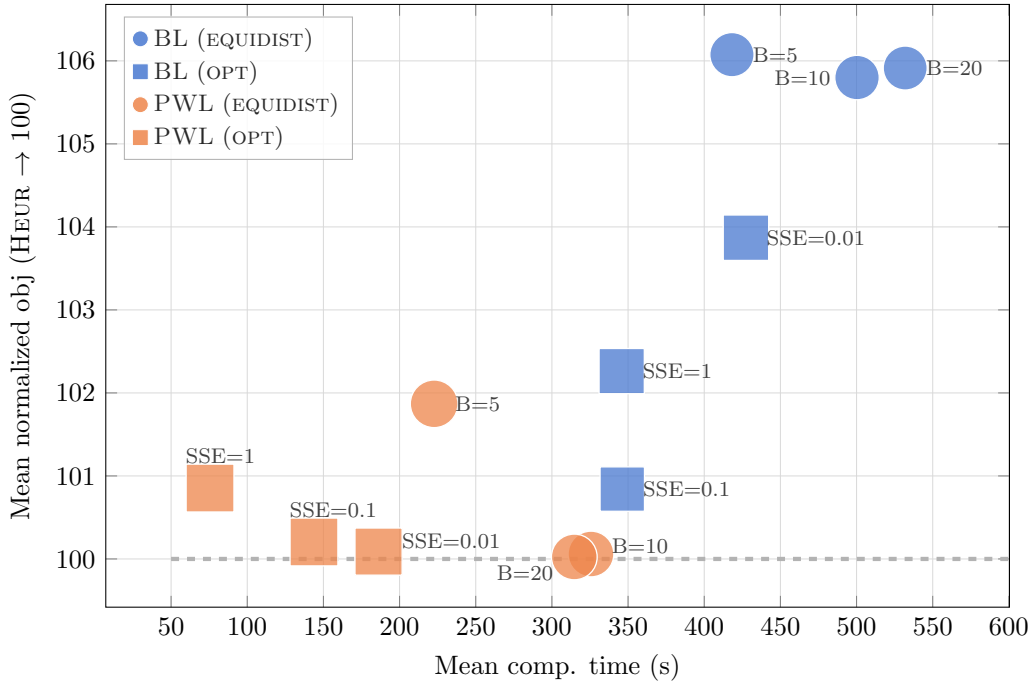
 
\autoref{fig:tradeoff_all_qos_baseline} shows a clear trade-off between computational effort and solution quality. As expected, increasing breakpoint resolution increases solve times. The effect on solution quality, however, differs between \pwlBLmodel and \pwlmodel. 
For \pwlmodel models, improvements in the objective values vary modestly, indicating that even rough approximations capture the relevant structure. 
In contrast, \pwlBLmodel is considerably more sensitive to the approximation setting, with recovered objectives ranging from approximately 0.8\% to 6.1\% above the baseline. 

Counterintuitively, the figure shows that for the \pwlBLmodel models, the more refined approximations do not always lead to better solutions than the coarser approximations. This is because all instances are included in the analysis, including those solved with a substantial optimality gap. When restricting the analysis to instances solved within 1\% of optimality, this effect is as expected; more refined models yield better solutions, ranging from 0.2\% above baseline for \pwlBLmodel using $B = 20$ to 3.9\% above baseline for \pwlBLmodel with $B = 5$. However, this comes at the cost of reduced robustness, as only 16 instances can be solved to near-optimality within the time limit in the setting yielding the best solutions. This analysis is available in \autoref{apx:approx}. 

Within \pwlmodel, the best configuration depends on the desired trade-off between speed, solution quality, and robustness. The optimized breakpoint strategy with SSE=1 is the most efficient, finding solutions to 61 instances and achieving the lowest average solve time. 
Increasing the number of equidistant breakpoints to $B = 20$ yields the best objective (0.02\% above the baseline), but finds solutions to only 52 instances within the time limit. The \opt strategy with SSE = 0.01 offers the most balanced approach, achieving an average objective value 0.1\% above baseline while solving 58 instances. 
Overall, the \equidist strategy is generally dominated by the \opt strategy, highlighting the value of adaptive breakpoint placement.

\section{Conclusion} \label{sec:conclusion}

In this paper, we study the \gls{SAMCF} problem in which demand responds endogenously to service quality while routing decisions are subject to hard capacity constraints. We formulate the problem as a nonlinear program that captures the interaction between routing, demand, and service quality, and propose a column-generation-based fixed-point heuristic to solve it. 

The computational study demonstrates that the proposed heuristic is highly effective. Across all test instances, it consistently finds high-quality solutions in a fraction of a second, significantly outperforming the nonlinear formulation and its piecewise-linear approximations in terms of robustness and scalability. In particular, while these models often struggle to solve large instances within the time limit, the heuristic provides near-optimal solutions and reliable performance across all demand profiles. These results highlight the practical value of combining column generation with iterative demand adjustment in service-aware network models.

From a modeling perspective, the results underscore the importance of capturing endogenous demand in capacitated network design problems. Allowing demand to respond to service quality fundamentally changes the allocation of capacity across commodities and can lead to substantially different routing decisions compared to classical inelastic formulations.
Several directions for future research arise from this work. First, it would be interesting to compare the proposed framework with alternative approaches based on generating sets of paths for each commodity a priori and subsequently solving the routing problem. Second, extending the model to incorporate congestion-dependent costs would allow for a more realistic representation of service quality and introduce an additional form of demand elasticity, as demand responds not only to routing decisions but also to congestion effects induced by aggregate flows. In this context, allowing for commodity-specific sensitivity parameters would enable heterogeneous responses to both service quality and congestion. Third, the short computation time of the proposed heuristic enables the \gls{SAMCF} framework to be integrated within network design as a tool to evaluate candidate network configurations. Fourth, the proposed methodology could be applied in other domains, including logistics and telecommunications, where demand and service levels are tightly coupled. Finally, in a transportation setting, the service-aware layer of the framework allows network operators to identify which demand is cost-effective to serve, while implicitly revealing demand segments that may be better suited for more flexible, e.g., on-demand, services.

    \printbibliography

    \newpage
    \appendix
    \renewcommand{\thesection}{Appendix~\Alph{section}}
    \section{Comparison of Routing Assumptions: Threshold vs SAMCF}\label{apx:EJORobjectives}
This section compares the routing assumptions used in \textcite{hansenExactAlgorithmPublic2026} with those adopted in this study. 
The instances considered are solutions generated as part of a public transport network design study \parencite{hansenExactAlgorithmPublic2026}, where passenger assignment relies on an OD-specific threshold: passengers are assumed to travel in the public transport network if the provided service meets a predefined quality threshold. 
The routing in \textcite{hansenExactAlgorithmPublic2026} relied on a set of pre-generated paths, so we reproduce it here using the column-generation approach described for solving the MCF problem in \autoref{sec:methods}. The total routing costs are very similar, differing by less than 0.4\%, indicating that the pre-generated path set is sufficient for these instances. However, this approach is not easily scalable. 

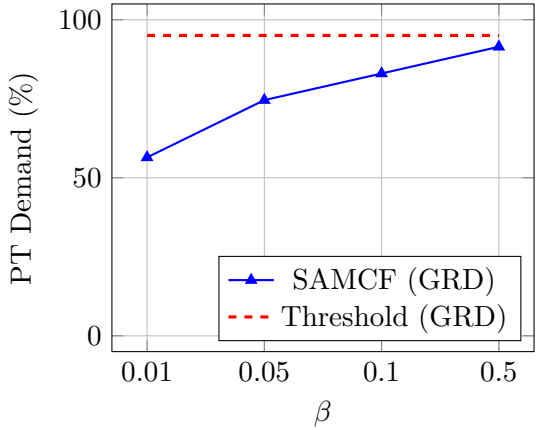
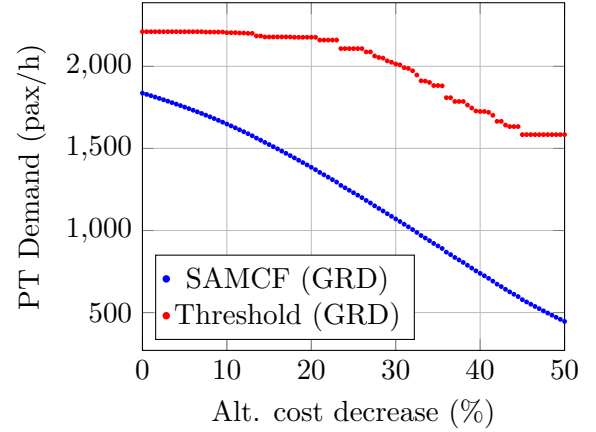
\begin{figure}[ht]
\centering
\begin{subfigure}{0.45\textwidth}
\centering
\begin{tikzpicture}
\begin{axis}[
    width=\textwidth,
    xlabel={$\beta$},
    ylabel={PT Demand (\%)},
    grid=major,
    ymin=0,
    ymax=100,
    enlarge y limits=0.05,
    legend pos=south east,
    legend columns=1,
    symbolic x coords={0.01,0.05,0.1,0.5},
    xtick=data,
]
\addplot[thick, mark=triangle*, blue] table[col sep=comma, x=Beta, y=DE_Beta_QoS_PctDemand_PT]
{de_beta_sweep_pt_demand_and_objective_grid.csv};
\addlegendentry{\gls{SAMCF} (\TPCodeFromInstance{Grid})}
\addplot[
    dashed,
    very thick,
    red
] table[col sep=comma, x=Beta, y=DE_Inelastic_PctDemand_PT_Baseline]
{de_beta_sweep_pt_demand_and_objective_grid.csv};
\addlegendentry{Threshold (\TPCodeFromInstance{Grid})}
\end{axis}
\end{tikzpicture}
\caption{Demand routed in PT network vs the logit-parameter $\beta$}
\label{apx:ejorA}
\end{subfigure}%
\hfill
\begin{subfigure}{0.45\textwidth}
\centering
\begin{tikzpicture}
\begin{axis}[
    width=\textwidth,
    xlabel={Alt. cost decrease (\%)},
    ylabel={PT Demand (pax/h)},
    grid=major,
    xmin=0,
    xmax=50,
    legend pos=south west,
    legend image post style={mark size=1pt}
]
\addplot[
    thick,
    mark=*,
    mark size=0.5pt,
    blue,
    only marks
] table[
    col sep=comma,
    x=AltModeCostDecreasePct,
    y=QoS_PT_Demand_Abs
] {pt_demand_abs_vs_altmode_cost_decrease_grid.csv};
\addlegendentry{\gls{SAMCF} (\TPCodeFromInstance{Grid})}
\addplot[
    thick,
    mark=*,
    mark size=0.5pt,
    red,
    only marks
] table[
    col sep=comma,
    x=AltModeCostDecreasePct,
    y=Inelastic_PT_Demand_Abs
] {pt_demand_abs_vs_altmode_cost_decrease_grid.csv};
\addlegendentry{Threshold (\TPCodeFromInstance{Grid})}
\end{axis}
\end{tikzpicture}
\caption{Demand routed in PT network vs the cost of the alternative mode}
\label{apx:ejorB}
\end{subfigure}
\caption{Comparison of the threshold-based model from \textcite{hansenExactAlgorithmPublic2026} to the demand response in this work, considering a target set of 300 OD pairs.}
\end{figure}

The networks were constructed so that a target subset of 300 OD pairs could be served in terms of network capacity, but not all OD pairs necessarily had a path with a cost below the threshold. As a result, not all demand is routed under this assumption.  
\autoref{apx:ejorA} illustrates the resulting share of demand routed in the public transport (PT) network as a function of the sensitivity parameter $\beta$, using a line plan generated for the \TPCodeFromInstance{Grid} instance. Under the threshold-based assumption, OD pairs switch entirely between modes once the threshold is met, resulting in a high PT demand. In contrast, the \gls{SAMCF} routing used in this study (computed via the proposed heuristic) routes less demand, as mode choice is governed by a continuous demand function rather than an all-or-nothing threshold. 

\autoref{apx:ejorB} shows the relationship between routed PT demand and the cost of the alternative mode, which provides further insight into the differences between the two routing approaches. As the alternative mode becomes cheaper relative to PT, demand shifts away from PT in both cases. In the threshold-based model, this occurs in discrete steps, with groups of OD pairs switching modes at once, whereas in the \gls{SAMCF} model, the transition is more gradual.

Although the results are shown for a single instance generated for the \TPCodeFromInstance{Grid} network, we observe the same general pattern across the other networks as well.

Overall, this comparison shows that the routing assumption has a notable impact on the passenger assignment and illustrates how the threshold-based approach produces an all-or-nothing switching behavior, while the \gls{SAMCF} model leads to a more gradual response in demand to changes in service. 

\section{Detailed Instance-Level Results for $\beta\in \{0.01, 0.1, 0.5\}$} \label{apx:detailed_tables}

The tables in this section report detailed results for each instance and model variant, including the objective value, computation time, optimality gap where applicable, and demand routed in the PT network. Since the \pwlmodel is an approximation of \nlmodel, we report the recovered objective value and demand routed in the PT network to allow for comparison across methods. Note that the gap reported for the \pwlmodel refers to the optimality gap reported by Gurobi and not the gap between the lower bound and the recovered objective value. 

\begin{table}[htbp]
\caption{Instance-level comparison of the heuristic (\heurShort), \pwlmodel, and \nlmodel for $\beta = 0.01$.}
\resizebox{\textwidth}{!}{%
\begin{tabular}{@{}r 
            @{\hspace{8mm}} rrr
            @{\hspace{8mm}}
            rrr
            @{\hspace{8mm}}
            rr
            @{\hspace{8mm}}
            rrr@{}}
\toprule
  & \multicolumn{3}{c}{Comp. time (s)} & \multicolumn{3}{c}{Objective value} & \multicolumn{2}{c}{Lower bound} & \multicolumn{3}{c}{PT demand} \\ \cmidrule(lr){2-4} \cmidrule(lr){5-7} \cmidrule(lr){8-9} \cmidrule(lr){10-12}   
  \#ODs & \heurShort & \nlmodel & \pwlmodel & \heurShort & \nlmodel & \pwlmodel & \nlmodel & \pwlmodel & \heurShort & \nlmodel & \pwlmodel \\ 
 \addlinespace[3pt]
 \multicolumn{12}{l}{\textbf{\TPCodeFromInstance{Athens}}} \\
172 & 0.19 & (7.63\%) & 8.15 & 73,670.90 & 74,535.04 & 73,726.67 & 68,846.38 & 72,087.08 & 647.40 & 601.39 & 642.93 \\
453 & 0.25 & -- & 56.14 & 121,112.91 & -- & 121,197.81 & 112,429.09 & 118,575.14 & 958.30 & -- & 950.00 \\
982 & 0.54 & -- & 387.34 & 172,413.72 & -- & 172,536.77 & 159,053.21 & 168,842.29 & 1,270.70 & -- & 1,256.31 \\
2,385 & 1.16 & -- & -- & 232,520.27 & -- & -- & 214,305.86 & -- & 1,585.00 & -- & -- \\
\addlinespace[3pt]
 \multicolumn{12}{l}{\textbf{\TPCodeFromInstance{Grid}}} \\
39 & 0.15 & (5.42\%) & 2.78 & 82,295.59 & 82,295.57 & 82,295.58 & 77,833.42 & 80,523.66 & 581.31 & 581.31 & 581.31 \\
101 & 0.06 & (6.84\%) & 34.29 & 120,568.17 & 120,568.13 & 120,568.15 & 112,326.91 & 117,952.85 & 863.17 & 863.17 & 863.17 \\
228 & 0.1 & (7.40\%) & 134.49 & 168,899.55 & 168,899.50 & 168,899.53 & 156,403.89 & 165,261.46 & 1,151.74 & 1,151.74 & 1,151.74 \\
567 & 0.3 & -- & (0.34\%) & 229,210.10 & -- & 229,227.51 & 211,944.35 & 223,670.12 & 1,436.13 & -- & 1,435.52 \\
\addlinespace[3pt]
 \multicolumn{12}{l}{\textbf{\TPCodeFromInstance{Lowersaxony}}} \\
12 & 0.05 & 5.54 & 0.02 & 28,698.57 & 28,698.57 & 28,698.57 & 28,697.12 & 28,085.10 & 376.48 & 376.48 & 376.48 \\
31 & 0.05 & (1.48\%) & 0.34 & 45,400.04 & 45,400.04 & 45,400.04 & 44,729.10 & 44,468.92 & 562.34 & 562.34 & 562.34 \\
78 & 0.06 & (4.52\%) & 2.34 & 66,565.36 & 66,565.35 & 66,565.35 & 63,558.21 & 65,211.79 & 754.66 & 754.66 & 754.66 \\
395 & 0.26 & -- & 29.13 & 93,596.73 & -- & 93,596.72 & 88,307.82 & 91,643.18 & 951.46 & -- & 951.46 \\
\addlinespace[3pt]
 \multicolumn{12}{l}{\textbf{\TPCodeFromInstance{Sioux}}} \\
61 & 0.06 & (6.10\%) & 0.61 & 100,676.58 & 100,676.57 & 100,676.57 & 94,531.67 & 98,308.24 & 937.61 & 937.61 & 937.61 \\
116 & 0.09 & (7.68\%) & 10.59 & 168,592.73 & 168,592.70 & 168,592.70 & 155,638.23 & 164,759.55 & 1,412.50 & 1,412.50 & 1,412.50 \\
209 & 0.11 & (8.38\%) & 81.41 & 240,532.87 & 240,532.83 & 240,532.83 & 220,384.19 & 235,116.52 & 1,892.09 & 1,892.09 & 1,892.09 \\
552 & 0.33 & -- & (0.17\%) & 321,821.72 & -- & 321,821.67 & 294,273.67 & 314,104.89 & 2,369.32 & -- & 2,369.32 \\
\bottomrule
\end{tabular}
}
\label{tab:beta0p01}
\end{table}

\begin{table}[htbp]
\caption{Instance-level comparison of the heuristic (\heurShort), \pwlmodel, and \nlmodel for $\beta = 0.1$.}
\resizebox{\textwidth}{!}{%
\begin{tabular}{@{}r 
            @{\hspace{8mm}} rrr
            @{\hspace{8mm}}
            rrr
            @{\hspace{8mm}}
            rr
            @{\hspace{8mm}}
            rrr@{}}
\toprule
  & \multicolumn{3}{c}{Comp. time (s)} & \multicolumn{3}{c}{Objective value} & \multicolumn{2}{c}{Lower bound} & \multicolumn{3}{c}{PT demand} \\ \cmidrule(lr){2-4} \cmidrule(lr){5-7} \cmidrule(lr){8-9} \cmidrule(lr){10-12}   
  \#ODs & \heurShort & \nlmodel & \pwlmodel & \heurShort & \nlmodel & \pwlmodel & \nlmodel & \pwlmodel & \heurShort & \nlmodel & \pwlmodel \\ 
 \addlinespace[3pt]
 \multicolumn{12}{l}{\textbf{\TPCodeFromInstance{Athens}}} \\
172 & 0.1 & -- & 145.68 & 58,921.11 & -- & 58,996.39 & 53,562.80 & 56,981.16 & 889.54 & -- & 888.50 \\
453 & 0.24 & -- & (0.27\%) & 98,513.21 & -- & 98,650.92 & 90,663.17 & 95,455.48 & 1,317.24 & -- & 1,313.87 \\
982 & 0.46 & -- & (0.91\%) & 141,699.15 & -- & 142,526.56 & 131,624.76 & 137,817.62 & 1,750.81 & -- & 1,718.15 \\
2,385 & 1.22 & -- & -- & 193,418.84 & -- & -- & 181,425.77 & -- & 2,191.19 & -- & -- \\
\addlinespace[3pt]
 \multicolumn{12}{l}{\textbf{\TPCodeFromInstance{Grid}}} \\
39 & 0.05 & (4.53\%) & 9.98 & 66,651.06 & 66,651.04 & 66,727.50 & 63,630.85 & 64,849.66 & 858.33 & 858.33 & 857.30 \\
101 & 0.08 & (6.46\%) & 141.1 & 97,668.64 & 97,668.60 & 97,780.10 & 91,356.01 & 94,921.58 & 1,264.00 & 1,264.00 & 1,262.50 \\
228 & 0.09 & (14.60\%) & 390.1 & 138,069.43 & 151,037.23 & 138,211.54 & 128,979.68 & 134,360.98 & 1,687.38 & 1,562.94 & 1,685.49 \\
567 & 0.16 & -- & -- & 190,227.82 & -- & -- & 178,495.89 & -- & 2,110.60 & -- & -- \\
\addlinespace[3pt]
 \multicolumn{12}{l}{\textbf{\TPCodeFromInstance{Lowersaxony}}} \\
12 & 0.04 & 28.42 & 0.14 & 24,930.71 & 24,930.69 & 24,938.79 & 24,930.69 & 23,507.49 & 458.53 & 458.53 & 458.40 \\
31 & 0.05 & (9.66\%) & 2.41 & 40,127.55 & 40,127.53 & 40,147.59 & 36,249.52 & 37,648.17 & 654.95 & 654.95 & 654.65 \\
78 & 0.15 & (11.62\%) & 15 & 57,640.04 & 57,640.01 & 57,678.94 & 50,939.64 & 54,449.05 & 903.53 & 903.53 & 902.98 \\
395 & 0.17 & (16.71\%) & (0.94\%) & 79,630.80 & 84,936.52 & 79,692.91 & 70,746.78 & 75,198.65 & 1,180.47 & 1,118.48 & 1,179.59 \\
\addlinespace[3pt]
 \multicolumn{12}{l}{\textbf{\TPCodeFromInstance{Sioux}}} \\
61 & 0.05 & (7.51\%) & 14.01 & 75,569.85 & 75,569.81 & 75,677.80 & 69,891.21 & 73,548.58 & 1,393.26 & 1,393.26 & 1,391.75 \\
116 & 0.07 & (6.97\%) & 37.72 & 128,096.89 & 128,096.85 & 128,291.48 & 119,168.50 & 125,282.39 & 2,116.99 & 2,116.99 & 2,114.37 \\
209 & 0.09 & (6.38\%) & 180.01 & 184,138.61 & 184,138.55 & 184,421.71 & 172,382.52 & 180,635.14 & 2,848.49 & 2,848.49 & 2,844.64 \\
552 & 0.28 & -- & (0.04\%) & 249,707.32 & -- & 250,077.37 & 235,409.89 & 245,294.19 & 3,577.75 & -- & 3,572.72 \\
\bottomrule
\end{tabular}
}
\label{tab:beta0p1}
\end{table}

\begin{table}[htbp]
\caption{Instance-level comparison of the heuristic (\heurShort), \pwlmodel, and \nlmodel for $\beta = 0.5$.}
\resizebox{\textwidth}{!}{%
\begin{tabular}{@{}r 
            @{\hspace{8mm}} rrr
            @{\hspace{8mm}}
            rrr
            @{\hspace{8mm}}
            rr
            @{\hspace{8mm}}
            rrr@{}}
\toprule
  & \multicolumn{3}{c}{Comp. time (s)} & \multicolumn{3}{c}{Objective value} & \multicolumn{2}{c}{Lower bound} & \multicolumn{3}{c}{PT demand} \\ \cmidrule(lr){2-4} \cmidrule(lr){5-7} \cmidrule(lr){8-9} \cmidrule(lr){10-12}   
  \#ODs & \heurShort & \nlmodel & \pwlmodel & \heurShort & \nlmodel & \pwlmodel & \nlmodel & \pwlmodel & \heurShort & \nlmodel & \pwlmodel \\ 
 \addlinespace[3pt]
\multicolumn{12}{l}{\textbf{\TPCodeFromInstance{Athens}}} \\
172 & 0.15 & -- & 89.94 & 56,054.93 & -- & 56,170.15 & 52,792.82 & 53,832.27 & 943.11 & -- & 940.74 \\
453 & 0.27 & -- & (0.02\%) & 94,641.20 & -- & 94,807.75 & 89,682.74 & 91,347.28 & 1,382.14 & -- & 1,378.84 \\
982 & 0.48 & -- & (0.08\%) & 137,101.44 & -- & 137,325.77 & 131,012.10 & 132,863.25 & 1,825.82 & -- & 1,824.77 \\
2,385 & 0.97 & -- & -- & 188,051.53 & -- & -- & 180,409.82 & -- & 2,275.98 & -- & -- \\
\addlinespace[3pt]
 \multicolumn{12}{l}{\textbf{\TPCodeFromInstance{Grid}}} \\
39 & 0.04 & (1.15\%) & 0.39 & 63,205.69 & 63,205.63 & 63,318.88 & 62,479.00 & 62,107.61 & 954.54 & 954.55 & 952.33 \\
101 & 0.06 & (2.58\%) & 1.57 & 92,754.82 & 92,754.72 & 92,915.27 & 90,357.84 & 90,870.78 & 1,404.76 & 1,404.76 & 1,401.68 \\
228 & 0.07 & -- & 3.64 & 131,881.44 & -- & 132,089.82 & 128,476.07 & 129,187.33 & 1,863.14 & -- & 1,859.22 \\
567 & 0.21 & -- & 21.09 & 182,905.82 & -- & 183,165.10 & 178,402.83 & 179,493.93 & 2,307.77 & -- & 2,302.92 \\
\addlinespace[3pt]
 \multicolumn{12}{l}{\textbf{\TPCodeFromInstance{Lowersaxony}}} \\
12 & 0.04 & 33.42 & 0.45 & 22,192.78 & 22,192.74 & 22,222.92 & 22,192.13 & 20,198.25 & 536.98 & 536.98 & 536.33 \\
31 & 0.05 & (10.34\%) & 6.09 & 37,104.26 & 37,104.19 & 37,144.33 & 33,266.01 & 33,449.26 & 721.20 & 721.20 & 720.38 \\
78 & 0.06 & -- & 38.55 & 53,792.89 & -- & 53,859.77 & 47,768.75 & 49,276.53 & 982.74 & -- & 981.38 \\
395 & 0.28 & -- & (0.47\%) & 75,069.02 & -- & 75,170.63 & 66,866.65 & 69,635.91 & 1,268.81 & -- & 1,266.81 \\
\addlinespace[3pt]
 \multicolumn{12}{l}{\textbf{\TPCodeFromInstance{Sioux}}} \\
61 & 0.05 & 223.79 & 0.47 & 69,463.80 & 69,463.76 & 69,655.29 & 69,461.70 & 68,301.00 & 1,572.19 & 1,572.19 & 1,568.29 \\
116 & 0.05 & (0.95\%) & 1.19 & 119,806.26 & 119,806.20 & 120,095.94 & 118,669.82 & 118,307.25 & 2,351.89 & 2,351.90 & 2,346.24 \\
209 & 0.07 & (19.49\%) & 6.7 & 173,842.31 & 213,035.06 & 174,233.53 & 171,513.80 & 172,047.31 & 3,133.51 & 2,112.15 & 3,126.05 \\
552 & 0.16 & -- & 13.22 & 237,679.00 & -- & 238,165.17 & 234,228.33 & 235,289.09 & 3,902.97 & -- & 3,893.89 \\
\bottomrule
\end{tabular}
}
\label{tab:beta0p5}
\end{table}

\section{Approximation Trade-Offs for the Subset of Models Solved to Near-Optimality} \label{apx:approx}

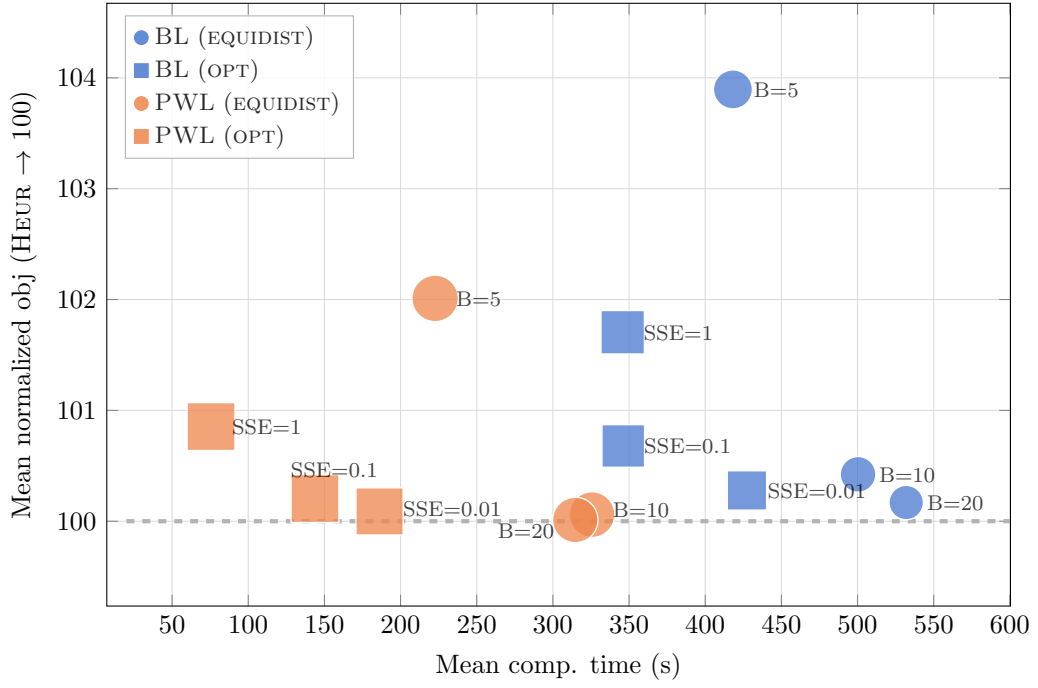
\begin{figure}[h]
\centering
\begin{tikzpicture}
\begin{axis}[
  width=0.85\linewidth,
  height=0.60\linewidth,
  xlabel={Mean comp. time (s)},
  ylabel={Mean normalized obj (\heurShort $\to$ 100)},
  grid=major,
  grid style={line width=0.3pt, draw=gray!30},
  tick label style={font=\small},
  label style={font=\small},
  legend style={at={(0.02,0.98)}, anchor=north west, font=\footnotesize,
    draw=gray!60, fill=white, fill opacity=0.85},
    legend cell align=left,
  clip=false,
  enlarge x limits=0.15,
  enlarge y limits=0.20,
]
\draw[dashed, gray!60, line width=1.6pt]
  (axis cs:20,100) -- (axis cs:600,100);

  \addplot[
    only marks, mark=*, mark size=7.32pt,
    mark options={fill=colBL, fill opacity=0.75, draw=white, line width=0.5pt},
    point meta=explicit symbolic,
    nodes near coords, nodes near coords style={font=\scriptsize, anchor=west, xshift=4pt},
    forget plot,
  ] coordinates {(418.1859, 103.8959) [B=5]};

  \addplot[
    only marks, mark=*, mark size=6.7pt,
    mark options={fill=colBL, fill opacity=0.75, draw=white, line width=0.5pt},
    point meta=explicit symbolic,
    nodes near coords, nodes near coords style={font=\scriptsize, anchor=west, xshift=4pt},
    forget plot,
  ] coordinates {(500.3012, 100.4229) [B=10]};

  \addplot[
    only marks, mark=*, mark size=6.46pt,
    mark options={fill=colBL, fill opacity=0.75, draw=white, line width=0.5pt},
    point meta=explicit symbolic,
    nodes near coords, nodes near coords style={font=\scriptsize, anchor=west, xshift=4pt},
    forget plot,
  ] coordinates {(531.9291, 100.1685) [B=20]};

  \addplot[
    only marks, mark=square*, mark size=7.39pt,
    mark options={fill=colBL, fill opacity=0.75, draw=white, line width=0.5pt},
    point meta=explicit symbolic,
    nodes near coords, nodes near coords style={font=\scriptsize, anchor=west, xshift=4pt},
    forget plot,
  ] coordinates {(427.3878, 100.2777) [SSE=0.01]};

  \addplot[
    only marks, mark=square*, mark size=8.01pt,
    mark options={fill=colBL, fill opacity=0.75, draw=white, line width=0.5pt},
    point meta=explicit symbolic,
    nodes near coords, nodes near coords style={font=\scriptsize, anchor=west, xshift=4pt},
    forget plot,
  ] coordinates {(346.1952, 100.6792) [SSE=0.1]};

  \addplot[
    only marks, mark=square*, mark size=8.17pt,
    mark options={fill=colBL, fill opacity=0.75, draw=white, line width=0.5pt},
    point meta=explicit symbolic,
    nodes near coords, nodes near coords style={font=\scriptsize, anchor=west, xshift=4pt},
    forget plot,
  ] coordinates {(345.8176, 101.7049) [SSE=1]};

  \addplot[
    only marks, mark=*, mark size=8.74pt,
    mark options={fill=colMcCormick, fill opacity=0.75, draw=white, line width=0.5pt},
    point meta=explicit symbolic,
    nodes near coords, nodes near coords style={font=\scriptsize, anchor=west, xshift=4pt},
    forget plot,
  ] coordinates {(222.6825, 102.0104) [B=5]};

  \addplot[
    only marks, mark=*, mark size=8.6pt,
    mark options={fill=colMcCormick, fill opacity=0.75, draw=white, line width=0.5pt},
    point meta=explicit symbolic,
    nodes near coords, nodes near coords style={font=\scriptsize, anchor=west, xshift=4pt,yshift=2pt},
    forget plot,
  ] coordinates {(325.6529, 100.0595) [B=10]};

  \addplot[
    only marks, mark=*, mark size=8.55pt,
    mark options={fill=colMcCormick, fill opacity=0.75, draw=white, line width=0.5pt},
    point meta=explicit symbolic,
    nodes near coords, nodes near coords style={font=\scriptsize, anchor=west, xshift=-33pt,yshift=-4pt},
    forget plot,
  ] coordinates {(314.7085, 100.0115) [B=20]};

  \addplot[
    only marks, mark=square*, mark size=8.87pt,
    mark options={fill=colMcCormick, fill opacity=0.75, draw=white, line width=0.5pt},
    point meta=explicit symbolic,
    nodes near coords, nodes near coords style={font=\scriptsize, anchor=west, xshift=5pt,yshift=1pt},
    forget plot,
  ] coordinates {(186.1988, 100.0887) [SSE=0.01]};

  \addplot[
    only marks, mark=square*, mark size=9.0pt,
    mark options={fill=colMcCormick, fill opacity=0.75, draw=white, line width=0.5pt},
    point meta=explicit symbolic,
    nodes near coords, nodes near coords style={font=\scriptsize, anchor=west, xshift=-13pt,yshift=11pt},
    forget plot,
  ] coordinates {(143.8114, 100.2044) [SSE=0.1]};

  \addplot[
    only marks, mark=square*, mark size=9.0pt,
    mark options={fill=colMcCormick, fill opacity=0.75, draw=white, line width=0.5pt},
    point meta=explicit symbolic,
    nodes near coords, nodes near coords style={font=\scriptsize, anchor=west, xshift=4pt},
    forget plot,
  ] coordinates {(75.5787, 100.8540) [SSE=1]};

\addlegendimage{only marks, mark=*, mark size=3pt,
  mark options={fill=colBL, draw=white}}
\addlegendentry{\pwlBLmodel (\equidist)}
\addlegendimage{only marks, mark=square*, mark size=3pt,
  mark options={fill=colBL, draw=white}}
\addlegendentry{\pwlBLmodel (\opt)}
\addlegendimage{only marks, mark=*, mark size=3pt,
  mark options={fill=colMcCormick, draw=white}}
\addlegendentry{\pwlmodel (\equidist)}
\addlegendimage{only marks, mark=square*, mark size=3pt,
  mark options={fill=colMcCormick, draw=white}}
\addlegendentry{\pwlmodel (\opt)}
\end{axis}
\end{tikzpicture}
\caption{Trade-off between computation time and solution quality for different breakpoint strategies considering the subset of instances solved within 1\% of optimality for each setting. The x-axis shows the average computation time in seconds, the y-axis shows the average objective value of the recovered, \nlmodel-feasible solutions obtained from the linear (square markers) and bilinear (round markers) approximations. 
The objective values have been normalized according to \heurShort $= 100$, and marker size indicates the number of instances solved within the optimality gap. }
\label{fig:tradeoff_solved1pct_qos_baseline}
\end{figure}

\autoref{fig:tradeoff_solved1pct_qos_baseline} illustrates the trade-off between computation time, solution quality, and robustness for the subset of instances solved within 1\% of optimality by the \pwlmodel or \pwlBLmodel.

\section{Numerical Results for Experiments with Outer-Bound PWL-Relaxations} \label{apx:twohour}

This section presents the numerical results for experiments using the outer PWL-relaxations. The \pwlBLmodel and \pwlmodel are formulated using an outer PWL approximation of the demand function, which guarantees an overestimation of the captured demand and therefore provides valid lower bounds for the original problem. These models, together with the \nlmodel, are evaluated with a time limit of up to two hours.

The resulting bounds are reported below. For reference, we include the solution obtained by the heuristic (\heurShort). For the relaxation models (\pwlBLmodel and \pwlmodel), we report the solver-provided objective value (``Relax. obj'') and lower bound (``LB''), as well as the objective value of the recovered solution that is feasible for the original nonlinear problem (``Recov. obj''). The reported gap corresponds to the optimality gap reported by the solver. 

\autoref{tab:2h_beta0p01_and_beta0p05} presents the results for $\beta \in \{0.01, 0.05\}$, while \autoref{tab:2h_beta0p1_and_beta0p5} presents the results for $\beta \in \{0.1, 0.5\}$.

\begin{sidewaystable}
\caption{Numerical results for experiments allowing up to two hours to solve the nonlinear model (\nlmodel), bilinear relaxation (\pwlBLmodel), and fully linearized relaxation (\pwlmodel)}
\label{tab:2h_beta0p01_and_beta0p05}
\centering
\resizebox{0.9\textwidth}{!}{%
\begin{tabular}{@{}lr@{\hspace{8mm}}rrr@{\hspace{8mm}}rrrr@{\hspace{8mm}}rrrr@{}}
\toprule
  & \heurShort & \multicolumn{3}{c}{\nlmodel} & \multicolumn{4}{c}{\pwlBLmodel} & \multicolumn{4}{c}{\pwlmodel} \\  \cmidrule(lr){2-2} \cmidrule(lr){3-5} \cmidrule(lr){6-9} \cmidrule(lr){10-13}   
 \#ODs & Obj & Obj & LB & Gap (\%) & Recov. obj & Relax. obj & LB & Gap (\%) & Recov. obj & Relax. obj & LB & Gap (\%) \\
\addlinespace[3pt]
\multicolumn{13}{l}{\textbf{\TPCodeFromInstance{Athens}} $\mathbf{(\beta = 0.01)}$} \\
172 & 73,670.90 & 73,670.89 & 69,133.10 & 6.16\% & 73,670.88 & 73,633.73 & 73,626.95 & 0.01\% & 73,726.67 & 72,207.28 & 72,207.28 & 0.00\% \\
 452 & 121,112.91 & 121,112.90 & 112,711.78 & 6.94\% & 127,092.69 & 127,058.32 & 120,997.60 & 4.77\% & 121,197.81 & 118,591.57 & 118,580.06 & 0.01\% \\
 982 & 172,413.72 & 173,989.19 & 159,754.82 & 8.18\% & 222,743.85 & 222,742.20 & 171,875.51 & 22.84\% & 172,536.77 & 168,780.86 & 168,764.82 & 0.01\% \\
 2,385 & 232,520.27 & 306,942.61 & 215,007.48 & 29.95\% &  &  & 231,719.04 &  &  &  &  &  \\
 \multicolumn{13}{l}{\textbf{\TPCodeFromInstance{Grid}} $\mathbf{(\beta = 0.01)}$} \\
 39 & 82,295.59 & 82,295.57 & 78,368.32 & 4.77\% & 82,295.57 & 82,256.27 & 82,248.65 & 0.01\% & 82,295.58 & 80,462.72 & 80,462.72 & 0.00\% \\
 101 & 120,568.17 & 120,568.13 & 112,823.04 & 6.42\% & 120,568.15 & 120,510.26 & 120,498.21 & 0.01\% & 120,568.15 & 117,951.31 & 117,951.31 & 0.00\% \\
 228 & 168,899.55 & 168,899.50 & 156,641.12 & 7.26\% & 168,899.52 & 168,818.76 & 168,799.43 & 0.01\% & 168,899.53 & 165,290.52 & 165,276.69 & 0.01\% \\
 567 & 229,210.10 & 229,905.15 & 211,950.43 & 7.81\% & 229,210.37 & 229,107.45 & 226,754.36 & 1.03\% & 229,210.07 & 224,386.77 & 224,068.18 & 0.14\% \\
 \multicolumn{13}{l}{\textbf{\TPCodeFromInstance{Lowersaxony}} $\mathbf{(\beta = 0.01)}$} \\
 12 & 28,698.57 & 28,698.57 & 28,697.10 & 0.01\% & 28,698.57 & 28,695.61 & 28,694.93 & 0.00\% & 28,698.57 & 28,133.49 & 28,133.49 & 0.00\% \\
 31 & 45,400.04 & 45,400.04 & 44,999.16 & 0.88\% & 45,400.04 & 45,392.36 & 45,388.82 & 0.01\% & 45,400.04 & 44,522.27 & 44,522.27 & 0.00\% \\
 78 & 66,565.36 & 66,565.35 & 63,894.82 & 4.01\% & 66,570.32 & 66,554.84 & 66,548.48 & 0.01\% & 66,565.35 & 65,202.70 & 65,200.96 & 0.00\% \\
 395 & 93,596.73 & 93,596.72 & 88,307.82 & 5.65\% & 93,751.07 & 93,725.48 & 93,499.63 & 0.24\% & 93,596.72 & 91,654.32 & 91,645.32 & 0.01\% \\
 \multicolumn{13}{l}{\textbf{\TPCodeFromInstance{Sioux}} $\mathbf{(\beta = 0.01)}$} \\
 61 & 100,676.58 & 100,676.57 & 95,064.04 & 5.57\% & 100,676.57 & 100,624.21 & 100,614.81 & 0.01\% & 100,676.57 & 98,541.72 & 98,541.72 & 0.00\% \\
 116 & 168,592.73 & 168,592.70 & 156,407.16 & 7.23\% & 168,592.70 & 168,504.61 & 168,488.30 & 0.01\% & 168,592.70 & 164,677.70 & 164,677.70 & 0.00\% \\
 209 & 240,532.87 & 240,532.83 & 221,011.79 & 8.12\% & 240,532.83 & 240,404.53 & 240,359.14 & 0.02\% & 240,532.83 & 234,727.89 & 234,727.89 & 0.00\% \\
 552 & 321,821.72 & 322,248.32 & 294,469.28 & 8.62\% & 324,429.25 & 324,266.57 & 319,967.61 & 1.33\% & 321,821.66 & 313,934.09 & 313,902.80 & 0.01\% \\
 \multicolumn{13}{l}{\textbf{\TPCodeFromInstance{Athens}} $\mathbf{(\beta = 0.05)}$} \\
172 & 63,064.39 & 63,064.37 & 55,615.13 & 11.81\% & 63,077.92 & 63,001.30 & 62,876.37 & 0.20\% & 63,064.37 & 60,616.28 & 60,610.52 & 0.01\% \\
453 & 104,676.39 & 104,680.87 & 93,106.54 & 11.06\% & 110,808.74 & 110,757.59 & 103,416.77 & 6.63\% & 104,699.15 & 100,940.73 & 100,935.73 & 0.00\% \\
982 & 149,979.92 & 231,754.24 & 134,545.75 & 41.94\% & 214,806.41 & 214,788.94 & 147,448.23 & 31.35\% & 150,027.36 & 144,894.58 & 144,799.61 & 0.07\% \\
2,385 & 203,752.71 &  & 184,411.72 &  &  &  & 200,083.14 &  &  &  &  &  \\
 \multicolumn{13}{l}{\textbf{\TPCodeFromInstance{Grid}} $\mathbf{(\beta = 0.05)}$} \\
39 & 71,071.39 & 71,071.37 & 66,638.78 & 6.24\% & 71,071.37 & 70,990.29 & 70,983.81 & 0.01\% & 71,071.37 & 68,408.30 & 68,407.68 & 0.00\% \\
101 & 104,027.46 & 104,027.42 & 94,696.58 & 8.97\% & 104,027.42 & 103,905.42 & 103,867.99 & 0.04\% & 104,027.42 & 100,175.81 & 100,166.53 & 0.01\% \\
228 & 146,586.59 & 146,586.54 & 132,506.29 & 9.61\% & 146,677.51 & 146,519.88 & 145,220.36 & 0.89\% & 146,586.54 & 141,263.41 & 141,250.15 & 0.01\% \\
567 & 201,139.20 & 206,762.04 & 182,582.76 & 11.69\% & 222,318.84 & 222,163.21 & 197,864.41 & 10.94\% &  &  &  &  \\
 \multicolumn{13}{l}{\textbf{\TPCodeFromInstance{Lowersaxony}} $\mathbf{(\beta = 0.05)}$} \\
12 & 26,471.49 & 26,471.47 & 26,471.47 & 0.00\% & 26,471.48 & 26,448.63 & 26,448.27 & 0.00\% & 26,471.48 & 25,207.75 & 25,207.75 & 0.00\% \\
31 & 42,006.21 & 42,006.19 & 41,211.68 & 1.89\% & 42,006.19 & 41,971.65 & 41,968.18 & 0.01\% & 42,006.20 & 39,988.93 & 39,988.93 & 0.00\% \\
78 & 60,511.29 & 60,511.27 & 56,028.53 & 7.41\% & 60,511.28 & 60,459.39 & 60,453.42 & 0.01\% & 60,511.27 & 57,756.75 & 57,756.75 & 0.00\% \\
395 & 83,929.54 & 83,929.50 & 75,590.57 & 9.94\% & 86,833.22 & 86,770.54 & 83,371.06 & 3.92\% & 83,932.26 & 80,357.43 & 79,701.80 & 0.82\% \\
 \multicolumn{13}{l}{\textbf{\TPCodeFromInstance{Sioux} $\mathbf{(\beta = 0.05)}$}} \\
61 & 82,887.36 & 82,887.33 & 74,283.73 & 10.38\% & 82,887.33 & 82,794.95 & 82,787.11 & 0.01\% & 82,887.34 & 79,698.26 & 79,698.26 & 0.00\% \\
116 & 139,181.67 & 139,181.63 & 124,216.01 & 10.75\% & 139,181.62 & 139,033.14 & 139,019.33 & 0.01\% & 139,181.63 & 134,271.72 & 134,265.33 & 0.00\% \\
209 & 199,087.56 & 199,087.49 & 177,795.39 & 10.69\% & 199,087.49 & 198,886.17 & 198,142.24 & 0.37\% & 199,087.50 & 192,468.84 & 192,450.06 & 0.01\% \\
552 & 268,568.11 & 268,956.35 & 240,672.20 & 10.52\% & 272,539.81 & 272,295.28 & 265,121.95 & 2.63\% & 268,568.03 & 259,970.37 & 259,911.76 & 0.02\% \\
\bottomrule
\end{tabular}}
\end{sidewaystable}

\begin{sidewaystable}
\caption{Numerical results for experiments allowing up to two hours to solve the nonlinear model (\nlmodel), bilinear relaxation (\pwlBLmodel), and fully linearized relaxation (\pwlmodel)}
\label{tab:2h_beta0p1_and_beta0p5}
\centering
\resizebox{0.9\textwidth}{!}{%
\begin{tabular}{@{}lr@{\hspace{8mm}}rrr@{\hspace{8mm}}rrrr@{\hspace{8mm}}rrrr@{}}
\toprule
  & \heurShort & \multicolumn{3}{c}{\nlmodel} & \multicolumn{4}{c}{\pwlBLmodel} & \multicolumn{4}{c}{\pwlmodel} \\  \cmidrule(lr){2-2} \cmidrule(lr){3-5} \cmidrule(lr){6-9} \cmidrule(lr){10-13}   
 \#ODs & Obj & Obj & LB & Gap (\%) & Recov. obj & Relax. obj & LB & Gap (\%) & Recov. obj & Relax. obj & LB & Gap (\%) \\
\addlinespace[3pt]
\multicolumn{13}{l}{\textbf{\TPCodeFromInstance{Athens}} $\mathbf{(\beta = 0.1)}$} \\
172 & 58,921.11 & 58,921.13 & 53,771.21 & 8.74\% & 58,965.99 & 58,871.24 & 58,810.95 & 0.10\% & 58,921.09 & 56,620.38 & 56,615.34 & 0.01\% \\
453 & 98,513.21 & 98,513.17 & 90,848.92 & 7.78\% & 98,964.55 & 98,822.12 & 97,428.32 & 1.41\% & 98,513.19 & 95,134.28 & 95,124.80 & 0.01\% \\
982 & 141,699.15 &  & 132,041.74 &  & 201,672.33 & 201,639.73 & 139,537.47 & 30.80\% & 141,758.24 & 137,328.47 & 137,278.01 & 0.04\% \\
2,385 & 193,418.84 &  & 181,800.14 &  &  &  & 190,705.24 &  &  &  &  &  \\
\multicolumn{13}{l}{\textbf{\TPCodeFromInstance{Grid}} $\mathbf{(\beta = 0.1)}$} \\
39 & 66,651.06 & 66,651.04 & 64,689.90 & 2.94\% & 66,651.04 & 66,543.67 & 66,538.42 & 0.01\% & 66,651.04 & 64,556.26 & 64,552.87 & 0.01\% \\
101 & 97,668.64 & 97,668.60 & 91,865.19 & 5.94\% & 97,668.60 & 97,511.86 & 97,502.12 & 0.01\% & 97,668.61 & 94,489.43 & 94,480.03 & 0.01\% \\
228 & 138,069.43 & 138,069.37 & 129,388.40 & 6.29\% & 138,069.37 & 137,854.96 & 137,208.33 & 0.47\% &  &  &  &  \\
567 & 190,227.82 &  & 179,757.61 &  & 199,639.10 & 199,400.03 & 187,632.76 & 5.90\% &  &  &  &  \\
\multicolumn{13}{l}{\textbf{\TPCodeFromInstance{Lowersaxony}} $\mathbf{(\beta = 0.1)}$} \\
12 & 24,930.71 & 24,930.69 & 24,930.69 & 0.00\% & 24,930.70 & 24,901.28 & 24,900.91 & 0.00\% & 24,930.70 & 22,938.01 & 22,938.01 & 0.00\% \\
31 & 40,127.55 & 40,127.53 & 36,875.40 & 8.10\% & 40,127.53 & 40,069.94 & 40,067.08 & 0.01\% & 40,127.53 & 37,016.31 & 37,013.04 & 0.01\% \\
78 & 57,640.04 & 57,640.01 & 51,882.63 & 9.99\% & 57,640.01 & 57,554.73 & 57,548.99 & 0.01\% & 57,640.02 & 53,711.80 & 53,707.53 & 0.01\% \\
395 & 79,630.80 & 79,637.59 & 70,343.84 & 11.67\% & 85,354.54 & 85,316.22 & 79,441.05 & 6.89\% & 79,630.77 & 75,072.39 & 74,895.36 & 0.24\% \\
\multicolumn{13}{l}{\textbf{\TPCodeFromInstance{Sioux}} $\mathbf{(\beta = 0.1)}$} \\
61 & 75,569.85 & 75,569.82 & 70,940.95 & 6.13\% & 75,569.82 & 75,441.50 & 75,438.24 & 0.00\% & 75,569.82 & 72,969.08 & 72,969.08 & 0.00\% \\
116 & 128,096.89 & 128,096.85 & 119,799.96 & 6.48\% & 128,098.31 & 127,879.37 & 127,842.28 & 0.03\% & 128,096.86 & 124,465.56 & 124,454.45 & 0.01\% \\
209 & 184,138.61 & 184,138.55 & 172,976.02 & 6.06\% & 184,265.78 & 183,954.13 & 183,728.23 & 0.12\% & 184,138.56 & 179,550.24 & 179,533.05 & 0.01\% \\
552 & 249,707.32 & 249,708.94 & 235,345.72 & 5.75\% & 288,834.77 & 288,523.26 & 248,009.22 & 14.04\% & 249,707.26 & 244,005.16 & 243,975.13 & 0.01\% \\
\multicolumn{13}{l}{\textbf{\TPCodeFromInstance{Athens}} $\mathbf{(\beta = 0.5)}$} \\
172 & 56,054.93 &  & 52,925.19 &  & 56,067.27 & 56,029.61 & 55,941.03 & 0.16\% & 56,054.92 & 53,655.80 & 53,653.20 & 0.00\% \\
453 & 94,641.20 &  & 89,931.40 &  & 94,795.38 & 94,739.70 & 94,312.59 & 0.45\% &  &  &  &  \\
982 & 137,101.44 &  & 130,998.23 &  & 137,713.51 & 137,676.58 & 136,473.49 & 0.87\% &  &  &  &  \\
2,385 & 188,051.53 &  & 180,725.24 &  &  &  &  &  &  &  &  &  \\
\multicolumn{13}{l}{\textbf{\TPCodeFromInstance{Grid}} $\mathbf{(\beta = 0.5)}$} \\
39 & 63,205.69 & 63,205.63 & 62,630.79 & 0.91\% & 63,205.63 & 63,161.76 & 63,157.64 & 0.01\% & 63,205.67 & 61,944.44 & 61,943.02 & 0.00\% \\
101 & 92,754.82 & 92,754.96 & 90,399.16 & 2.54\% & 92,756.38 & 92,701.98 & 92,695.78 & 0.01\% & 92,754.79 & 90,638.58 & 90,637.33 & 0.00\% \\
228 & 131,881.44 & 131,881.30 & 128,381.37 & 2.65\% & 131,881.32 & 131,814.83 & 131,420.01 & 0.30\% & 131,881.38 & 128,898.19 & 128,885.74 & 0.01\% \\
567 & 182,905.82 &  & 178,401.03 &  & 184,025.96 & 183,975.62 & 181,614.24 & 1.28\% & 182,905.75 & 179,146.42 & 179,128.69 & 0.01\% \\
\multicolumn{13}{l}{\textbf{\TPCodeFromInstance{Lowersaxony}} $\mathbf{(\beta = 0.5)}$} \\
12 & 22,192.78 & 22,192.74 & 22,192.74 & 0.00\% & 22,192.75 & 22,168.92 & 22,166.84 & 0.01\% & 22,192.77 & 20,037.25 & 20,037.25 & 0.00\% \\
31 & 37,104.26 & 37,104.19 & 34,611.28 & 6.72\% & 37,104.18 & 37,053.58 & 37,050.01 & 0.01\% & 37,104.24 & 33,202.44 & 33,200.95 & 0.00\% \\
78 & 53,792.89 & 53,792.79 & 48,645.44 & 9.57\% & 53,792.80 & 53,727.18 & 53,672.45 & 0.10\% & 53,792.85 & 48,967.04 & 48,964.36 & 0.01\% \\
395 & 75,069.02 & 84,814.65 & 67,577.27 & 20.32\% & 75,084.35 & 75,012.40 & 74,731.95 & 0.37\% & 75,068.99 & 69,581.43 & 69,462.24 & 0.17\% \\
\multicolumn{13}{l}{\textbf{\TPCodeFromInstance{Sioux}} $\mathbf{(\beta = 0.5)}$} \\
61 & 69,463.80 & 69,463.76 & 69,459.32 & 0.01\% & 69,463.77 & 69,410.15 & 69,403.43 & 0.01\% & 69,463.78 & 68,034.47 & 68,034.47 & 0.00\% \\
116 & 119,806.26 & 119,806.30 & 119,794.53 & 0.01\% & 119,806.21 & 119,735.62 & 119,726.56 & 0.01\% & 119,806.23 & 117,912.00 & 117,902.24 & 0.01\% \\
209 & 173,842.31 & 173,842.25 & 171,716.48 & 1.22\% & 173,842.25 & 173,757.73 & 173,745.27 & 0.01\% & 173,842.27 & 171,548.46 & 171,531.49 & 0.01\% \\
552 & 237,679.00 & 275,593.06 & 234,228.33 & 15.01\% & 237,678.92 & 237,576.78 & 237,544.44 & 0.01\% & 237,678.96 & 234,683.50 & 234,661.17 & 0.01\% \\
\bottomrule
\end{tabular}}
\end{sidewaystable}

\end{document}